\documentclass{gtart_h}
\def\ifplaintex{\expandafter\ifx\csname documentclass\endcsname\relax}

\def\gtp{{\mathsurround=0pt\it $\cal G\mskip-2mu$eometry \&\ 
$\cal T\!\!$opology $\cal P\!$ublications}}  % GT publications

\def\recd{{\small Received:\qua\receiveddate\ifx\reviseddate\relax
\else\qquad Revised:\qua\reviseddate\fi\par}} 

\def\lognumber#1{\def\thelognumber{#1}}
\def\volumenumber#1{\def\thevolumenumber{#1}}
\def\volumeyear#1{\def\thevolumeyear{#1}}
\def\papernumber#1{\def\thepapernumber{#1}}
\def\pagenumbers#1#2{\def\startpage{#1}\def\finishpage{#2}}
\def\published#1{\def\publishdate{#1}}

\def\received#1{\def\receiveddate{#1}}
\def\revised#1{\def\reviseddate{#1}}
\def\accepted#1{\def\accepteddate{#1}}
\def\asciititle#1{\def\theasciititle{#1}}
\def\covertitle#1{\def\thecovertitle{#1}}

\def\asciiaddress#1{\def\theasciiaddress{#1}}
\def\asciiemail#1{\def\theasciiemail{#1}}

\long\def\asciiabstract#1{\long\def\theasciiabstract{#1}}

\let\\\par\let\thelognumber\relax\let\thevolumenumber\relax
\let\thepapernumber\relax\let\thevolumeyear\relax\let\startpage\relax
\let\finishpage\relax\let\publishdate\relax\let\receiveddate\relax
\let\reviseddate\relax\let\accepteddate\relax\let\theasciititle\relax
\let\thecovertitle\relax\let\theasciiauthors\relax\let\theasciiaddress\relax
\let\theasciiabstract\relax

\let\theasciiemail\relax

\ifplaintex
\font\logobig=cmssbx10 scaled 3836
\font\logomed=cmssbx10 scaled 2557
\else
\font\logobig=cmssbx10 scaled 4200
\font\logomed=cmssbx10 scaled 2800
\fi

\long\def\makeagttitle{   %%% start of definition of \makeagttitle
\count0=\startpage
\agt\hfill      %   Journal title (top left) 
\hbox to 45truept{\vbox to 0pt{\vglue -13truept{\logomed A\kern -.37em{\logobig 
T}\kern -.38em G}\vss}\hss}
\break
{\small Volume \thevolumenumber\ (\thevolumeyear)
\startpage--\finishpage\nl
Published: \publishdate}

\vglue .25truein

{\parskip=0pt\leftskip 0pt plus
1fil\def\\{\par\smallskip}{\Large\bf\thetitle}\par\medskip} \vglue
0.05truein

{\parskip=0pt\leftskip 0pt plus 1fil\def\\{\par}{\sc\theauthors}
\par\medskip}%
 
\vglue 0.03truein 

{\small\leftskip 25truept\rightskip 25truept{\bf Abstract}\stdspace\theabstract

{\bf AMS Classification}\stdspace\theprimaryclass
\ifx\thesecondaryclass\relax\else; \thesecondaryclass\fi\par
{\bf Keywords}\stdspace \thekeywords\par}\vglue 7truept

}   %%%% end of definition of \makeagttitle

\ifplaintex
\font\phead=cmsl9 scaled 950
\font\pnum=cmbx10 scaled 913
\font\pfoot=cmsl9 scaled 950
\headline{\vbox to 0pt{\vskip -4.5mm\line{\small\phead\ifnum
\count0=\startpage ISSN 1472-2739 (on-line) 1472-2747 (printed)
\hfill {\pnum\folio}\else\ifodd\count0\def\\{ }% 
\ifx\theshorttitle\relax\thetitle\else\theshorttitle\fi\hfill{\pnum\folio}
\else\def\\{ and }{\pnum\folio}\hfill\ifx\theshortauthors\relax\theauthors
\else\theshortauthors\fi\fi\fi}\vss}}
\footline{\vbox to 0pt{\vglue 0mm\line{\small\pfoot\ifnum\count0=\startpage
\copyright\ \gtp\hfill\else
\agt, Volume \thevolumenumber\ (\thevolumeyear)\hfill\fi}\vss}}
\else
\font=cmsl9 scaled 1050
\font\lnum=cmbx10 
\font=cmsl9 scaled 1050
\makeatletter
\def\@oddhead{{\small\lhead\ifnum\count0=\startpage ISSN 1472-2739 
(on-line) 1472-2747 (printed)\hfill {\lnum\number\count0}\else\ifodd\count0
\def\\{ }\ifx\theshorttitle\relax \thetitle \else\theshorttitle\fi\hfill
{\lnum\number\count0}\else\def\\{ and }{\lnum\number\count0}
\hfill\ifx\theshortauthors\relax 
\theauthors\else\theshortauthors\fi\fi\fi}}\def\@evenhead{\@oddhead}
\def\@oddfoot{\small\lfoot\ifnum\count0=\startpage\copyright\ \gtp\hfill\else
\agt, Volume \thevolumenumber\ (\thevolumeyear)\hfill\fi}
\def\@evenfoot{\@oddfoot}
\makeatother
\fi
\let\maketitlepage\makeagttitle
\let\makeshorttitle\maketitlepage
\let\maketitle\maketitlepage

\newwrite\gtoutfile
\long\gdef\makeheadfile{  %%% start of definition of \makeheadfile
{\def\\{, }\def\s{ }
\immediate\openout\gtoutfile head.xxx
\immediate\write\gtoutfile{Proxy-for: \ifx\theasciiauthors\relax
\theauthors\else\theasciiauthors\fi\s<\ifx\theasciiemail\relax\theemail\else\theasciiemail\fi>}
\immediate\write\gtoutfile{\noexpand\\}
\immediate\write\gtoutfile{Authors: \ifx\theasciiauthors\relax
\theauthors\else\theasciiauthors\fi}
{\def\\{ }\immediate\write\gtoutfile{Title: \ifx\theasciititle\relax
\thetitle\else\theasciititle\fi}}
\immediate\write\gtoutfile{Subj-class: GT or SG, GR etc}
\immediate\write\gtoutfile{MSC-class: \theprimaryclass\ifx\thesecondaryclass\relax\else, \thesecondaryclass\fi}
\immediate\write\gtoutfile{Journal-ref: Algebr. Geom. Topol. \thevolumenumber\s
(\thevolumeyear) \startpage-\finishpage}
\immediate\write\gtoutfile{Comments: Published by Algebraic and
Geometric Topology at}
\immediate\write\gtoutfile{\s\s\s  http://www.maths.warwick.ac.uk/agt/AGTVol\thevolumenumber/agt-\thevolumenumber-\thepapernumber.abs.html}
\immediate\write\gtoutfile{\noexpand\\}
\immediate\write\gtoutfile{}
\ifx\theasciiabstract\relax
\immediate\write\gtoutfile{\theabstract}\else
\immediate\write\gtoutfile{\theasciiabstract}\fi
\immediate\write\gtoutfile{}
\immediate\write\gtoutfile{\noexpand\\}
\immediate\write\gtoutfile{}
\immediate\closeout\gtoutfile}}  %%% end of definition of \makeheadfile

\def\maketitlepage{\makeagttitle\makeheadfile}
\let\makeshorttitle\maketitlepage
\let\maketitle\maketitlepage

\lognumber{40}
\volumenumber{5}
\volumeyear{2005}
\papernumber{40}
\pagenumbers{965}{997}
\received{16 November 2004} 
\revised{30 May 2005}
\accepted{25 July 2005}
\published{11 August 2005}

\usepackage{amssymb,amsmath,amscd} 

\newtheorem{thm}{Theorem}[section]
\newtheorem{prop}[thm]{Proposition}
\newtheorem{lemma}[thm]{Lemma}
\newtheorem{cor}[thm]{Corollary}
\theoremstyle{definition}
    \newtheorem{definition}[thm]{Definition}
       \newtheorem{remark}[thm]{Remark}
       \newtheorem{example}[thm]{Example}
       \newtheorem{question}[thm]{Question}

\newcommand{\ol}[1]{\overline{#1} }

\newcommand{\E}{\boldsymbol{E}}

\newcommand{\bfa}{\mathbf{a}}
\newcommand{\bfA}{\mathbf{A}}

\newcommand{\bfx}{\mathbf{x}}

\newcommand{\xvec}[2]{#1_1,\ldots,#1_{#2}}

\newcommand{\cupd}{\overset{.}{\cup}}
\newcommand{\cupb}{\overset{b}{\cup}}

\newcommand{\R}{\ensuremath{R} }
\newcommand{\pres}[2]{\langle  #1  \vert  #2 \rangle}

\newcommand{\p}[1]{\pi_1(#1)}

\newcommand{\Span}{\operatorname{span}}
\newcommand{\Ad}{\operatorname{Ad}}
\newcommand{\image}{\operatorname{Im}}

\newcommand{\Rk}{\operatorname{rk}}

\newcommand{\tr}{\operatorname{tr}}
\newcommand{\Ker}{\operatorname{Ker}}

\newcommand{\tors}{\operatorname{tors}}

\renewcommand{\hom}{\operatorname{Hom}}

\newcommand{\TaZ}[2]{T^{\mathrm{Zar}}_{#1}(#2)}

\newcommand{\CC}{\mathbb{C}}

\newcommand{\HH}{\mathbb{H}}
\newcommand{\NN}{\mathbb{N}}

\newcommand{\RR}{\mathbb{R}}
\newcommand{\ZZ}{\mathbb{Z}}

\newcommand{\SU}[1][2]{\mathrm{SU} (#1)}
\newcommand{\U}[1][1]{\mathrm{U} (#1)}
\newcommand{\su}[1][2]{\mathfrak{su} (#1)}

\newcommand{\PSL}{\mathrm{PSL}_2}
\newcommand{\SL}{\mathrm{SL}_2}
\newcommand{\Sl}{\mathfrak{sl}_2}
\newcommand{\g}{\mathfrak{g}}

\newcommand{\ord}{\operatorname{ord}}

{\bigl(
\begin{smallmatrix} }%
{
\end{smallmatrix} \bigr)}

\begin{document}

\title{Deformations of reducible representations\\of 
3-manifold groups into $\mathrm{PSL}_{2}(\mathbb{C})$}
\covertitle{Deformations of reducible representations\\ 
of 3-manifold groups into PSL$_{2}(\noexpand\bf C)$}
\asciititle{Deformations of reducible representations
of 3-manifold groups into PSL_2(C)}
\shorttitle{Deformations of reducible representations}

\author{Michael Heusener\\Joan Porti}                  
\address{Laboratoire de Math\'ematiques, Universit\'e Blaise
Pascal\\63177, Aubi\`ere Cedex, France}                  
\secondaddress{Departament de Matem\`atiques, Universitat Aut\`onoma de
Barcelona\\E-08193 Bellaterra, Spain}

\asciiaddress{Laboratoire de Mathematiques, Universite Blaise
Pascal\\63177, Aubiere Cedex, France\\and\\Departament de 
Matematiques, Universitat Autonoma de
Barcelona\\E-08193 Bellaterra, Spain} 

\asciiemail{heusener@math.univ-bpclermont.fr, porti@mat.uab.es}           
\gtemail{\mailto{heusener@math.univ-bpclermont.fr}{\rm\qua
and\qua}\mailto{porti@mat.uab.es}}

\begin{abstract}
    Let $M$ be a 3-manifold with torus boundary which is a rational
    homology circle.
    We study deformations of reducible
     representations of $\p M$ into
    $\mathrm{PSL}_{2}(\mathbb{C})$ associated to a simple zero of the
    twisted Alexander polynomial.
    We also describe the local structure of the representation and
character varieties.
\end{abstract}

\asciiabstract{%
Let M a 3-manifold with torus boundary which is a rational homology
circle.  We study deformations of reducible representations of p_1(M)
into PSL_2(C) associated to a simple zero of the twisted Alexander
polynomial.  We also describe the local structure of the
representation and character varieties.}


\maketitle 

\cl{\small\sl Dedicated to the memory of Heiner Zieschang}
 
\section{Introduction}
\label{sec:intro}

This article is a continuation of the work started in \cite{HPS01def}.
Let $M$ be a connected, compact, orientable,  3-manifold
such that
$
\partial M$ is a torus. We assume that the first Betti number
$\beta_1(M)$ is one, i.e.\ $M$ is a rational homology circle.
In particular, $M$ is the exterior of a knot in a rational homology sphere.

Given a homomorphism $\alpha\co\p M\to\CC^*$, we define an abelian representation
$\rho_\alpha\co\p M\to \PSL(\CC)$ as follows:
\begin{equation}
\label{eq:rho-alpha}
\rho_\alpha(\gamma)=\pm
\begin{pmatrix} \alpha^{\frac12}(\gamma) &0\\ 0&
\alpha^{-\frac12}(\gamma)
\end{pmatrix}\qquad\forall\gamma\in\pi_1(M)
\end{equation}
where $\alpha^{\frac12}\co \p M\to\CC^{*}$ is a map (not
necessarily a homomorphism) such that
$(\alpha^{\frac12}(\gamma))^{2} =
\alpha(\gamma)$ for all $\gamma\in\p M$. The representation
$\rho_{\alpha}$ is reducible, i.e.\ $\rho_{\alpha}(\p M)$ has
global fixed points in $P^{1}(\CC)$.

\begin{question}
When can  $ \rho_{\alpha}$ be  deformed into irreducible
representations (i.e. representations whose images have no fixed points
in $P^{1}(\CC)$)?
\end{question}

Different versions of this question have been studied in
\cite{FK91}, \cite{Her97} and \cite{HK98} for $\SU $ and
\cite{BA98}, \cite{BA00}, \cite{BAL02},  \cite{HPS01def} and
\cite{Sho91} for $\SL (\CC)$.

The answer is related to a twisted Alexander invariant.
We first choose an isomorphism:
\begin{equation}
\label{eqn:explicitsplitting}
 H_1(M;\ZZ)\cong \tors(H_1(M;\ZZ))\oplus\ZZ,
\end{equation}
which amounts to choosing a projection onto the torsion subgroup
$H_1(M;\ZZ)\to \tors(H_1(M;\ZZ))$ and to fix a generator $\phi$ of
$H^1(M;\ZZ)\cong \hom(H_1(M;\ZZ),\ZZ)$. So, $\alpha$ induces a
homomorphism $\tors(H_1(M;\ZZ))\oplus\ZZ\to \CC^{*}$ which will be
also denoted by $\alpha$.

The composition of the
projection $\p M \to \tors(H_1(M;\ZZ))$ with the restriction of
$\alpha$ gives a representation
$\sigma\co \p M\to U(1)\subset\CC^*$:
\[
\sigma\co \p M\to \tors(H_1(M;\ZZ)) \overset\alpha\to  \CC^*.
\]
A homomorphism
$\phi_{\sigma}\co\pi_1(M)\to \CC[t^{\pm 1}]^*$
to the units of the ring of Laurent polynomials $\CC[t^{\pm 1}]$ is
given by
$\phi_{\sigma}(\gamma)=
\sigma(\gamma) t^{\phi(\gamma)}$.
 This allows a definition of the
 twisted Alexander polynomial $\Delta^{\phi_\sigma}(t)\in
 \CC[t^{\pm 1}]$, whose construction will be recalled in
Section~\ref{sec:Alex}.

We say that $\alpha$ is \emph{a zero of the Alexander invariant} if
  $\Delta^{\phi_\sigma}( \alpha(0,1) )=0$, where
$(0,1)\in \tors(H_1(M;\ZZ))\oplus\ZZ$. We show in
Section~\ref{sec:abelian} that being a zero and the  
 order of the
zero does not depend of the choice of the isomorphism
(\ref{eqn:explicitsplitting}).

We prove in Lemma~\ref{lem:nondef}
 that if $\rho_{\alpha}$ can be deformed into irreducible
representations, then $\alpha$ is a zero of the Alexander invariant.
For a simple zero this condition is also sufficient and we have
stronger conclusions, as the next theorem shows. Let $R(M)=\hom
(\pi_1(M),\PSL(\CC))$ denote the variety of representations of
$\pi_1(M)$ in $\PSL(\CC)$.

\begin{thm}
\label{thm:structureofR}
 If $\alpha$ is a simple zero of the Alexander invariant, then
$\rho_\alpha$ is contained in precisely two irreducible components
of $R(M)$, one of dimension 4 containing irreducible
representations and another of dimension 3 containing only abelian
ones. In addition, $\rho_\alpha$ is a smooth point of both
varieties and  the intersection  at the orbit of $\rho_\alpha$ is
transverse.
\end{thm}
When the representation $\alpha$ is trivial, then
it is not a zero of the Alexander invariant, because the Alexander
invariant is the usual untwisted Alexander polynomial $\Delta$,
and $\Delta(1)= \pm \vert \operatorname{tors} \,H_1(M,\mathbb
Z)\vert \neq 0$.

Let $X(M)=R(M)/\!/\PSL(\CC)$ denote the algebraic quotient where
$\PSL(\CC)$ acts by conjugation on $R(M)$.
The character of a representation $\rho\in R(M)$ is a map
$\chi_{\rho}\co\p M\to\CC$ given by $\chi_{\rho}(\gamma) = \tr^{2}
\rho(\gamma)$ for $\gamma\in\p M$.
There is a one-to-one correspondence between $X(M)$ and the set of
characters. Hence we call $X(M)$ the  variety
of characters of $\pi_1(M)$ (see \cite{HP04} for details).

Let $\chi_{\alpha}$ be the character of
$\rho_{\alpha}$.
\begin{thm}
\label{thm:structureofX} If $\alpha$ is a simple zero of the
Alexander invariant, then $\chi_\alpha$ is contained in precisely
two irreducible components of $X(M)$, which are curves and are the
quotients of the components of $R(M)$ in
Theorem~\ref{thm:structureofR}. In addition $\chi_\alpha$ is a
smooth point of both curves and the intersection at $\chi_\alpha$
is transverse.
\end{thm}

This paper generalizes the main results of
\cite{HPS01def} where we considered
only representations $\alpha\co\p M \to\CC^{*}$ which factor through
$H_{1}(M;\ZZ) / \tors(H_{1}(M;\ZZ))$ and for which
$\alpha^{\frac12}$ can be chosen as a homomorphism.
These conditions imply that $\rho_{\alpha}$ and its deformations can
be lifted to representations into $\SL(\CC)$. On the other hand, it was
shown in
\cite[Theorem~1.4]{HP04} that the representation variety $R(M)$ can
have many components which do not lift to the 
$\SL(\CC)$--representation variety.
 Hence Theorem \ref{thm:structureofR} and Theorem
\ref{thm:structureofX} generalize  the main results of
\cite{HPS01def}. Moreover, we have removed the condition that
$\rho_{\alpha}$ is not $\partial$-\emph{trivial} from \cite{HPS01def}.
Here a representation $\rho\in R(M)$ is called $\partial$-trivial if
\[
\rho\circ i_{\#}\co \p{\partial M}\to\PSL(\CC)
\]
is trivial.
An example will be given where the results of this paper apply but
those of \cite{HPS01def} do not.

In \cite{HPS01def} we considered  the usual Alexander polynomial,
but here we need a twisted version. The strategy of the
proof of Theorem~\ref{thm:structureofR} in this paper
is
similar to the one of Theorem~1.1 of \cite{HPS01def}: we construct
a metabelian representation $\rho^+\co\p M \to \PSL(\CC)$ which is
not abelian and has the same character as $\rho_{\alpha}$, and we
show that $\rho^+$ is a smooth point of $R(M)$. This involves
quite elaborate cohomology computations. More precisely, due to an
observation by Andr\'e Weil, the Zariski tangent space
$\TaZ{\rho}{R(M)}$ of $R(M)$ at a representation $\rho\in R(M)$
may  be viewed as a subspace of
the space of group cocycles $Z^1(\pi_1(M),\Sl(\CC)_{\rho})$. Here
$\Sl(\CC)_{\rho}$ denotes the $\p M$-module $\Sl(\CC)$ via
$\Ad\circ \rho$.
%is the tangent space to the orbit by conjugation.
The approach given here for these cohomological
computations and for the analysis of the tangent space is completely
self contained and simplifies in several aspects the computations
from \cite{HPS01def}. In particular, the new approach permits us to
remove the assumption that $\rho_\alpha$ is not $\partial$-trivial.

The paper is organized as follows. In Section~\ref{sec:Alex} we
recall the definition of the twisted Alexander polynomial and describe
its main properties. In Section~\ref{sec:notation} we recall some
basic facts from group cohomology and Weil's
construction which will be used in the sequel.
Section~\ref{sec:abelian} relates the vanishing of twisted
Alexander invariants to some elementary cohomology and
deformations of abelian representations. The next three sections
are devoted to prove that the metabelian representation $\rho^+$
can be deformed into irreducible representations. The cohomology
computations are done in Section~\ref{sec:cohomologymetabelian},
with a key lemma proved in Section~\ref{sec:cohomology}. The
smoothness of $R(M)$ at $\rho^+$  is proved in
Section~\ref{section:deformingmetabelian}.
Theorems~\ref{thm:structureofR} and \ref{thm:structureofX} are
proved in Sections~\ref{section:quadraticcone} and
\ref{section:characters} respectively. Finally,
Section~\ref{section:realvalued} is devoted to describe the local
structure of the set of real valued characters.

\rk{Acknowledgement} The authors would like to thank the referee for
his careful and thorough reading of the manuscript and for providing
useful suggestions for improving the paper.  The second author is
partially supported by the Spanish MCYT through grant BFM2003-03458.

\section{Twisted Alexander polynomial}
\label{sec:Alex}

Let $M$ be a manifold as in the introduction. We fix a projection
$p\co H_{1}(M;\ZZ)\to\tors(H_{1}(M;\ZZ))$ and a generator 
\[
\phi\in
H^{1}(M;\ZZ)=\hom(H_{1}(M;\ZZ),\ZZ) =\hom(\pi_{1}(M),\ZZ)\,,
\]
i.e.\ we fix an isomorphism as in (\ref{eqn:explicitsplitting})
\begin{displaymath}
\begin{array}{rcr}
H_{1}(M;\ZZ)&\cong&\tors(H_{1}(M;\ZZ))\oplus \ZZ\\
z&\mapsto &(p(z), \phi(z))\,.
\end{array}
\end{displaymath}
 
For every
representation $\sigma\co \tors(H_{1}(M;\ZZ))\to\U  \subset\CC^*$
the composition
$$ \p M  \to H_{1}(M;\ZZ) \stackrel{p}\longrightarrow
\tors(H_{1}(M;\ZZ))\stackrel{\sigma}\longrightarrow \U$$
will be  denoted by $\sigma(p) \co\p M\to \U$.
We consider the induced representation
\[
\begin{array}{rcl}
\phi_{\sigma(p)}\co\p M & \to & \CC[t^{\pm 1}]^* \\
\gamma&\mapsto & \sigma(p)(\gamma) t^{\phi(\gamma)}\,.
\end{array}
\]
In this way $\CC[t^{\pm 1}]$  is a $\pi_1(M)$-module (or a
$H_1(M;\ZZ)$-module since $\image \phi_{\sigma(p)}$ is abelian).

In the sequel we shall fix a projection 
$p\co H_{1}(M;\ZZ)\to\tors(H_{1}(M;\ZZ))$.
We shall, by convenient abuse of notation,  
continue to write $\sigma$ for $\sigma(p)$. 
Let 
\[
H_*^{\phi_{\sigma}}(M)\quad\textrm{ and }\quad
H^*_{\phi_{\sigma}}(M)
\]
denote the homology
and cohomology twisted by $\phi_{\sigma}$.
Using singular chains, $H_*^{\phi_{\sigma}}(M)$ and
$H^*_{\phi_{\sigma}}(M)$  can be defined respectively as the homology
and cohomology of the chain and cochain complexes:
$$
\CC[t^{\pm 1}]\otimes_{\pi_1(M)} C_*(\widetilde M;\ZZ)
\quad\textrm{ and } \quad \hom_{\pi_1(M)}(C_*(\widetilde
M;\ZZ),\CC[t^{\pm 1}]),
$$
where $\widetilde M$ denotes the universal covering of $M$. Alternatively, since
 $\phi_\sigma$ is abelian,
 we could take the maximal abelian covering of $M$ instead of the universal one.

In the sequel we shall write $R := \CC[t^{\pm 1}]$. Since $R$ is a
principal ideal domain and since $H_1^{\phi_{\sigma}}(M)$ is finitely
generated, we have a canonical decomposition
$$
H_1^{\phi_{\sigma}}(M)=R/r_0 R\oplus\cdots\oplus R/r_m R,
$$
where $r_i\in R$ and $r_{i+1}\mid r_i$.

\begin{definition}
The $R$-module $H_1^{\phi_{\sigma}}(M)$ is called the \emph{Alexander
  module} and 
$$
\Delta^{\phi_\sigma}_k= r_k\,r_{k+1}\,\cdots r_m
$$
the \emph{$k$-th twisted Alexander polynomial},  for
$k=0,\ldots,m$. The first one is also called the twisted Alexander
polynomial: $\Delta^{\phi_\sigma} :=\Delta^{\phi_\sigma}_0$.
\end{definition}
It is well defined up to units
in $R=\CC[t^{\pm 1}]$, i.e. up to
multiplication with elements $a\,t^n$ with $a\in \CC^*$ and
$n\in\ZZ$. We use the natural extension $\Delta^{\phi_\sigma}_k=1$
for $k>m$. 
Note that if $A\in M_{m,n}(R)$ is a presentation matrix for
$H_1^{\phi_{\sigma}}(M)$ then $\Delta^{\phi_\sigma}_k\in R$ is the
greatest common divisor of the minors of $A$ of order $(n-k)$.
%The definition of the Alexander polynomial and of the
Alexander module can be done in a more general context using only an
U.F.D. (see \cite[IV.3]{Tur03}).

\paragraph{Changing the isomorphism}
Note that $H_1^{\phi_{\sigma}}(M)$ and hence
$\Delta^{\phi_\sigma}_k$ are not invariants of the pair
$(M,\sigma)$, they depend on the  isomorphism
(\ref{eqn:explicitsplitting}); equivalently, 
they depend on the choice of the projection 
$p\co H_1(M,\ZZ)\to \tors(H_1(M;\ZZ))$ and the generator
$\phi\co H_1(M,\ZZ)\to\ZZ$.

Let $p_1,p_2\co H_1(M;\ZZ)\to \tors(H_1(M;\ZZ))$ be two projections.
They differ by a morphism 
$\psi\co \ZZ\to \tors(H_1(M;\ZZ))$.
Namely, for all $z\in H_1(M;\ZZ)$,
\begin{equation}
\label{eqn:p2=p1-psi}
p_2(z)= p_1(z)+\psi(\phi(z)).
\end{equation}
Therefore, given $\sigma\co\tors(H_1(M;Z)\to\U$, the
induced representations 
$\sigma_i := \sigma(p_i)\co\p M\to\U$ satisfy: 
\begin{equation}
\label{eqn:sigmaothersplitting}
\sigma_2 (\gamma)=\sigma_1 (\gamma)
\sigma(\psi(\phi(\gamma)))\qquad\forall\gamma\in\pi_1(M)\, .
\end{equation}
Hence, 
\[\phi_{\sigma_2}(\gamma) = \sigma_2(\gamma) t^{\phi(\gamma)} =
\sigma_1(\gamma)\sigma(\psi(\phi(\gamma)))t^{\phi(\gamma)}=
\sigma_1(\gamma)(a\, t)^{\phi(\gamma)}\] 
and $\phi_{\sigma_1}(\gamma)= \sigma_1(\gamma) t^{\phi(\gamma)}$
differ by replacing $t$ by 
$a\, t$, where $a=\sigma(\psi(1))\in\U$ and
$1$ denotes the
generator of $\ZZ$.
Therefore
\begin{equation}
\label{eqn:deltaothersplitting}
\Delta^{\phi_{\sigma_2}}_k(t) = 
\Delta^{\phi_{\sigma_1}}_k( a\, t)\,.
\end{equation}
The generator $\phi\co H_1(M;\ZZ)\to\ZZ$ is unique up to sign,
and 
replacing $\phi$  by $-\phi$ implies replacing $t$ by
$t^{-1}$ in the twisted polynomial.

\paragraph{Symmetry}
Consider the following involution on $\CC[t^{\pm 1}]$:
$$
\overline{\sum_i a_i\, t^i}=\sum_i \overline{a_i} t^{-i},
$$
where $\overline{a_i}$ denotes the complex conjugate of $a_i\in\CC$.
An ideal $I\subset\CC[t^{\pm 1}]$ is called \emph{symmetric} if
$I=\bar I$ and an element $\eta\in\CC[t^{\pm 1}]$ is called
\emph{symmetric} if it generates a symmetric ideal. Hence an element
$\eta\in\CC[t^{\pm 1}]$ is symmetric if and only if there exists a unit
$\epsilon\in\CC[t^{\pm 1}]^{*}$ such that $\bar\eta=\epsilon \eta $.
Notice that some authors use the expression \emph{weakly symmetric}
\cite{Tur03}.

\begin{prop}\label{pro:symmetry}
Let $M$ be a 3-manifold such that $\beta_{1}(M)=1$ and that $
\partial M$ is a torus.
For each homomorphism $\sigma\co \tors(H_{1}(M;\ZZ))\to \U$ and for
each splitting of (\ref{eqn:explicitsplitting}) we have that
$\Delta_{k}^{\phi_\sigma}$ is symmetric
i.e.,  $\Delta_{k}^{\phi_\sigma}$ and $\ol{\Delta_{k}^{\phi_\sigma}}$ are
equal up to multiplication by a unit of $\CC[t^{\pm 1}]$.
\end{prop}
\begin{proof}
Given a $R=\CC[t^{\pm 1}]$-module $N$,  $\ol N$ denotes the
$R$-module with the opposite $R$-action, i.e.\ $r \bar n := \bar r
\bar n$, for $r\in R$ and $\bar n\in \ol N$.
 Using the Blanchfield
duality pairing we obtain an isomorphism of $R$-modules
$D\co H^{\phi_{\sigma}}_{3-p}(M,
\partial M)\to
\ol{H^{p}_{\phi_{\sigma}}(M)}$. Since $R$ is a P.I.D., we obtain
\begin{gather}
\label{eq:duality}
 D\co \tors (H^{\phi_{\sigma}}_{n-p}(M,
\partial M) )\cong
         \ol{\tors (H_{p-1}^{\phi_{\sigma}}(M))} \\
\intertext{and }
\Rk_{R} H^{\phi_{\sigma}}_{n-p}(M,
\partial M)
                  = \Rk_{R} H^{\phi_{\sigma}}_{p}(M)\notag
\end{gather}
(see \cite{Bla57}, \cite{Mil62}, \cite[Sec.~2]{Lev77}
and \cite[Sect.~7]{Gor77}).

The proposition follows from the duality
formula~(\ref{eq:duality}) (see the proof of Theorem~7.7.1 in
\cite[p.~97]{Kaw96} for the details).
\end{proof}

\begin{remark} \label{rem:AlexNul}
In contrast to the untwisted situation, the Alexander module
$H_{1}^{\phi_\sigma}(M)$ can have nonzero rank. Examples are
easily obtained as follows: let $M_{1}$ be the complement of  knot
in  a homology sphere and let $M_{2}$ be a rational homology
sphere. Then $\p {M_{1}\# M_{2}}\cong\p {M_{1}}\ast\p{M_{2}}$ and
$ H_{1}(M_{1}\# M_{2};\ZZ)\cong H_{1}(M_{1};\ZZ) \oplus
H_{1}(M_{2};\ZZ)$ comes with a canonical splitting.

Since $H_{1}(M_{1};\ZZ)$ is torsion free, we can choose
$\phi\co \p {M_{1}\# M_{2}}\to\ZZ$ to be the composition 
\[ \phi\co \p {M_{1}\# M_{2}}\to H_{1}(M_{1}\# M_{2};\ZZ)\to
H_{1}(M_{1};\ZZ) \cong \ZZ\,.\]
For each nontrivial representation 
$\sigma\co H_{1}(M_{2};\ZZ)\to\U$ we obtain
$\phi_\sigma (h_{1}+ h_{2}) = \sigma(h_{2})
t^{\phi(h_{1})}\in\CC[t^{\pm 1}]^{*}$ for $h_{i}\in
H_{1}(M_{i};\ZZ)$ and hence 
$\sigma=\phi_\sigma\mid_{H_{1}(M_{2};\ZZ)}$. 

Since $\sigma$ is
nontrivial it follows that $H_{0}^{\phi_\sigma}(M_{2})=
H_{0}^{\phi_\sigma}(M_{1}\# M_{2}) = 0$. Moreover, we obtain that
$H_{0}^{\phi_\sigma}(S^{2})\cong\CC[t^{\pm 1}]$ and hence the
Mayer-Vietoris sequence gives a short exact sequence:
 \begin{equation*}
 0\to H_{1}^{\phi_\sigma}(M_{1})\oplus H_{1}^{\phi_\sigma}(M_{2})
  \to H_{1}^{\phi_\sigma}(M_{1}\# M_{2})\to \Ker j \to 0
\end{equation*}
where  $j\co H_{0}^{\phi_\sigma}(S^{2})\to
H_{0}^{\phi_\sigma}(M_{1})$
is surjective.
Since $H_{0}^{\phi_\sigma}(M_{1})$ is torsion it follows that
$\Ker j$ is a free $ \CC[t^{\pm 1}]$-module of rank one
and hence, $H_{1}^{\phi_\sigma}(M_{1}\# M_{2})$ has nonzero rank.
\end{remark}

\section{Group cohomology: Fox calculus and products}
\label{sec:notation}

Fox calculus will be used to compute the twisted Alexander
polynomial. Since 
this is a tool in group cohomology, we first
need the following lemma, that will also be used later. Details and
conventions 
in group cohomology can be found
in \cite{Bro82} and \cite{Wei94}.

\begin{lemma}\label{lem:EilMac}
Let $A$ be a $\p M$-module and let $X$ be any CW-complex with
$\pi_1(X)\cong \pi_1(M)$.
Then there are natural morphisms $H_i(X;A)\to H_i(\p M ;A)$ which are
isomorphisms
for $i=0,1$ and a surjection for $i=2$. In cohomology
there are natural morphisms $H^i(\p M;A)\to H^i( X ;A)$ which are
isomorphisms
for $i=0,1$ and an injection for $i=2$.
\end{lemma}
\begin{proof}
It is possible to construct an
Eilenberg-MacLane space $K$ of type $(\p M,1)$ from $X$ by attaching
$k$-cells, $k\geq 3$.
In this way we obtain a CW-pair $(K,X)$ and it follows that
$H_{j}(K,X;A)= 0$ and $H^{j}(K,X;A)= 0$  for  $j=1,2$ (this a direct
application of 
Theorems (4.4) and (4.4*) of
\cite[VI.4]{Whi78}). Hence
the exact sequences of the pair $(K,X)$ give the result.
\end{proof}

\paragraph{Fox calculus}

Let $\sigma\co\tors(H_1(M;\ZZ)\to\U$ be a representation and fix a
splitting 
of $H_1(M;\ZZ)$ as in (\ref{eqn:explicitsplitting}).
We can use Fox calculus to compute $\Delta_{k}^{\phi_\sigma}$:
choose a cell decomposition of $M$ with only one zero 
cell $x_{0}$. Since every presentation of $\p M$ obtained from a
cell decomposition of $M$ has deficiency one, we have:
 $$
\p M=\langle S_{1},\ldots,S_{n}\mid
R_{1},\ldots,R_{n-1}\rangle\,.$$
Denote by $\pi\co F_{n}\to\p M$ the canonical
projection where $F_{n}=F(S_{1},\ldots,S_{n})$ is the free group
generated by $n$ elements and
by $\partial/\partial
S_{i}\co \ZZ F_{n}\to\ZZ F_{n}$ the partial derivations of the group
ring of the free group.

The \emph{Jacobian}
$J^{\phi_{\sigma}}:= (J^{\phi_{\sigma}}_{ji})\in M_{n-1,n}(\CC[t^{\pm 1}])$ is defined
by  $J^{\phi_{\sigma}}_{ji}:=\phi_\sigma\circ\pi(
\partial R_{j}/
\partial S_{i})\in\CC[t^{\pm 1}]$.
Analogous to  \cite[Chapter~9]{BuZi03},
one can show that $J^{\phi_{\sigma}}$ is a presentation matrix for
$H^{\phi_\sigma}_{1}(M,x_{0})$.
The exact sequence for the pair
$(M,x_{0})$
yields $H^{\phi_\sigma}_{1}(M,x_{0})\cong
H^{\phi_\sigma}_{1}(M)\oplus\CC[t^{\pm 1}]$ (see \cite[pp.\ 61-62]{Tur03}).
Hence,
$\Delta_{k}^{\phi_\sigma}(M) =\Delta^{\phi_{\sigma}}_{k+1}(M,x_{0})$
and
$\Delta_{k}^{\phi_\sigma}(M)$ is the greatest common divisor of the
$(n-k-1)$-minors of the Jacobian $J^{\phi_{\sigma}}\in M_{n-1,n}(\CC[t^{\pm 1}])
$.
%\marnote{J'ai remplac\'e $\Delta^{\phi_{\sigma}}_{k}(M,x_{0})$ par
%$\Delta^{\phi_{\sigma}}_{k+1}(M,x_{0})$ (translation d'indice
%oblige).}

\begin{example}\label{ex:1erepartie}
 Let $M$ be the punctured torus bundle over
$S^1$ whose action of the
monodromy on $H_1(\dot{T}^2,\mathbb Z)$ is given by the
matrix
$\big(\begin{smallmatrix} 1 & 2 \\ 2 & 5
\end{smallmatrix} \big)$.
The fundamental group $\p{\dot T^2}$ is a free
group of rank $2$ generated by $\alpha$ and $\beta$.
A presentation of $\p M \cong \ZZ \ltimes \p{\dot T^2}$ is given by
\[
\p M=\langle \alpha,\beta,\mu\mid \mu\alpha\mu^{-1}=\alpha\beta^{2},
\mu\beta\mu^{-1}=\beta(\alpha\beta^{2})^{2}\rangle\,.
\]
Moreover, $H_{1}(M;Z)\cong (\ZZ/2 \oplus\ZZ/2)\oplus \ZZ$ comes with a
canonical splitting $s(1)=\mu$ i.e., $p(\mu)= 0$,
$p(\alpha)=\alpha$ and $p(\beta)=\beta$.
A generator $\phi\in H^{1}(M;\ZZ)\cong
\hom(\p
M;\ZZ)$ is given by $\phi(\mu) =1$ and $\phi(\alpha)=\phi(\beta)=0$.

There are exactly four
representations $\sigma_{i}\co(\ZZ/2 \oplus\ZZ/2)\to \U$
which give
rise to the following homomorphisms
$\phi_{\sigma_{i}}\co\p M\to\CC[t^{\pm 1}]^{*}$:
\begin{align*}
\phi_{\sigma_{1}}\co& \left\{
\begin{array}{lc}
\mu &\mapsto t\\
\alpha &\mapsto  1\\
\beta &\mapsto 1
\end{array}\right.  \, , \quad &
\phi_{\sigma_{2}}\co& \left\{
\begin{array}{lc}
\mu &\mapsto t\\
\alpha &\mapsto  1\\
\beta &\mapsto -1
\end{array}\right.  \, ,\\
\phi_{\sigma_{3}}\co& \left\{
\begin{array}{lc}
\mu &\mapsto t\\
\alpha &\mapsto  -1\\
\beta &\mapsto 1
\end{array}\right.  \, , \quad &
\phi_{\sigma_{4}}\co& \left\{
\begin{array}{lc}
\mu &\mapsto t\\
\alpha &\mapsto  -1\\
\beta &\mapsto -1 
\end{array}\right. \,.
\end{align*}
Now, a direct calculation gives:
\begin{align*}
J^{\phi_{\sigma_1}}&=
\begin{pmatrix}
0 & t-1 & -2\\
0 & -2 & t-5
\end{pmatrix},\quad&
J^{\phi_{\sigma_2}}&=
\begin{pmatrix}
0 & t-1 & 0\\
2 & 2 & t-1
\end{pmatrix}\\
J^{\phi_{\sigma_3}}&=
\begin{pmatrix}
2 & t-1 & 2\\
0 & 0 & t-1
\end{pmatrix},\quad &
J^{\phi_{\sigma_4}}&=
\begin{pmatrix}
2 & t-1 & 0\\
2 & 0 & t-1
\end{pmatrix}\,.
\end{align*}
Hence,
$\Delta^{\phi_{\sigma_{1}}}=t^{2}-6t+1$ and
$\Delta^{\phi_{\sigma_{i}}}=t-1$ for $i=2,3,4$.
\end{example}

\paragraph{Products in cohomology}
Let $\Gamma$ be a group and let $A$ be a $\Gamma$-module.
We denote by $(C^{*}(\Gamma;A),d)$ the
{normalized} cochain complex.
The coboundaries (respectively cocycles, cohomology) of $\Gamma$
with coefficients in $A$ are denoted by $B^{*}(\Gamma;A)$
(respectively $Z^{*}(\Gamma;A)$, $H^{*}(\Gamma;A)$).

Let $A_{1}$, $A_{2}$ and $A_{3}$ be $\Gamma$-modules.
The cup product of two cochains $u \in C^p(\Gamma;A_{1})$ and $v\in
C^q(\Gamma;A_{2})$ is the cochain
$u \cup v \in C^{p+q}(\Gamma;A_1\otimes A_2)$ defined by
\begin{equation}\label{eq:cupprod}
u \cup v(\xvec{\gamma}{p+q}) := u(\xvec{\gamma}{p}) \otimes
    \gamma_1\cdots\gamma_p
    \circ v(\gamma_{p+1},\ldots,\gamma_{p+q})\,.
\end{equation}
 Here $A_{1}\otimes A_{2}$ is a $\Gamma$-module via the diagonal
action.

It is possible to combine the cup product with any bilinear map
(compatible with the $\Gamma$ action) $b \co A_1\otimes A_2 \to A_3$.
So we obtain a cup product
$$
\cupb \co
     C^p(\Gamma;A_{1}) \otimes C^q(\Gamma;A_{2}) \to
                  C^{p+q}(\Gamma;A_3)\,. $$
                  For details see \cite[V.3]{Bro82}.
If $A=\g$ is a Lie algebra, then we obtain the \emph{cup-bracket}
of two cochains, which will be denoted by $[u \cup v]$. Note that the
cup-bracket is not associative on the cochain level. 
 If $A$ is an algebra then the cup product will be simply denoted by
 $\cupd$. This cup product 
is associative on the cochain level if the multiplication in $A$
is associative.

     Let $b\co A_{1}\otimes A_{2}\to A_{3}$ be bilinear  and
       let $z_i\in Z^1(\Gamma;A_{i})$, $i=1,2$, be cocycles. We define
       $f\co\Gamma\to A_{3}$ by $ f(\gamma):=
       b(z_1(\gamma)\otimes z_2(\gamma))$. A direct
       calculation gives: $d\,f(\gamma_1,\gamma_2)+
       b(z_1(\gamma_1)\otimes \gamma_1 \circ z_2(\gamma_2))
                + b(\gamma_1 \circ z_1(\gamma_2)\otimes
z_2(\gamma_1))=0$.
       This shows that:
        \begin{equation}\label{eq:cupsym}
              d\,f  + z_1 \cupb z_2 +
z_2 \overset{b\circ\tau}{\cup} z_1 =0
        \end{equation}
 where $\tau\co A_{1}\otimes A_{2}\to A_{2}\otimes A_{1}$ is the
twist operator.

\paragraph{Group cohomology and representation varieties}
Let $\Gamma$
be a group and let $\rho\co\Gamma\to\PSL(\CC)$ be a representation.
The Lie algebra $\Sl(\CC)$ turns into a $\Gamma$-module via
$\Ad\circ\rho$. We shall denote this $\Gamma$-module by
$\Sl(\CC)_{\rho}$. A cocycle $d\in Z^{1}(\Gamma;\Sl(\CC)_{\rho})$ is
a map $d\co\Gamma\to\Sl(\CC)_{\rho}$ satisfying
\[
d(\gamma_1\gamma_2)=d( \gamma_1)+\Ad_{\rho(\gamma_1)}d(\gamma_2),\;
\forall\gamma_1,\gamma_2\in\pi_1(M)\,.
\]

It was observed by Andr\'e Weil \cite{Wei64} that there is a
natural inclusion of the Zariski tangent space
$\TaZ{\rho}{R(\Gamma)}\hookrightarrow
Z^{1}(\Gamma;\Sl(\CC)_{\rho})$. Informally speaking, given a
smooth curve $\rho_{\epsilon}$ of representations through
$\rho_{0}=\rho$ one gets a $1$-cocycle
$d\co\Gamma\to\Sl(\CC)_{\rho}$ by defining
\[ d(\gamma) := \left.\frac{d \, \rho_{\epsilon}(\gamma)}
{d\,\epsilon}\right|_{\epsilon=0} \rho(\gamma)^{-1},
\quad\forall\gamma\in\Gamma\,.\]
It is easy to see that the
tangent space to the orbit by conjugation corresponds to the space
of $1$-coboundaries $B^{1}(\Gamma;\Sl(\CC)_{\rho})$ (see for
instance \cite[Section~4.5]{Kap01}). Here,
$b\co\Gamma\to\SL(\CC)_{\rho}$ is a coboundary if there exists
$x\in\Sl(\CC)$ such that $b(\gamma)=\Ad_{\rho(\gamma)}x-x$.
 A detailed account can be found in \cite[Thm.\ 2.6]{LM85} or
\cite[Ch.\ VI]{Rag72}. %%or \cite[\S 3.1.3]{Por97}.

 Let $\dim_{\rho} R(\Gamma)$ be the \emph{local dimension} of $R(\Gamma)$ at
$\rho$ (i.e. the maximal dimension of the irreducible components
of $R(\Gamma)$ containing $\rho$,
cf.\ \cite[Ch.\ II, \S1.4]{Sha77}). In the
sequel we shall use the following lemma:
\begin{lemma} \label{lem:locdim}
        Let $\rho\in\R(\Gamma)$ be given.
        If $\dim_{\rho} R(\Gamma)= \dim Z^{1}(\Gamma;\Sl(\CC)_{\rho})$
        then $\rho$ is a smooth
        point of the representation variety $\R(\Gamma)$ and $\rho$ is
        contained
        in a unique component of $R(\Gamma)$ of dimension $\dim
        Z^1(\Gamma;\Sl(\CC)_\rho)$.
\end{lemma}
\begin{proof}
        For every $\rho\in\R(\Gamma)$ we have
\[
\dim_{\rho} \R(\Gamma) \leq \dim \TaZ{\rho}{\R(\Gamma)}
        \leq \dim Z^1(\Gamma;\Sl(\CC)_{\rho})\,.
\]
The lemma follows from the fact that the equality
$\dim_{\rho} \R(\Gamma) = \dim \TaZ{\rho}{\R(\Gamma)}$ is the
condition in algebraic geometry  that guarantees that $\rho$
belongs to a single irreducible component of $R(\Gamma)$ and it is
a smooth point (for more details see \cite[Ch.\ II]{Sha77}).
\end{proof}

 \section{Abelian representations and the twisted Alexander invariant}
\label{sec:abelian}

Let $\alpha\co\p M\to\CC^{*}$ be a homomorphism. In what follows, the induced
homomorphism $H_{1}(M;Z)\to\CC^{*}$ will be also denoted by $\alpha$
and we denote by $\sigma$ the restriction of $\alpha$ to the
torsion subgroup, i.e.\ 
\[
\sigma := \alpha\mid_{\tors(H_1(M;\ZZ))}\co
\tors(H_1(M;\ZZ))\to\U\,.
\] 
Let us fix an isomorphism $ H_1(M;\ZZ)\cong
\tors(H_1(M;\ZZ))\oplus\ZZ$ as in 
(\ref{eqn:explicitsplitting}), i.e.\ we fix an projection 
$p\co H_1(M;\ZZ)\to\tors(H_1(M;\ZZ))$ and a generator $\phi\in
H^1(M;\ZZ)$.
The induced section $s_p\co \ZZ\to H_1(M;\ZZ)$ satisfies
\begin{displaymath}
s_p\circ\phi=\operatorname{Id}-p.
\end{displaymath}
%\marnote{J'ai remplac\'e $\phi_\sigma -p$ par
%$Id-p$.}

\begin{definition}
\label{definition:zero}
 We say that $\alpha$ is a \emph{zero of the $k$-th Alexander
   invariant of order $r$}
if $a :=\alpha(s_p(1))\in\CC^{*}$ is a zero of $\Delta_k^{\phi_{\sigma}}(t)$
of order $r$.
\end{definition}

\begin{lemma}
This definition does not depend on the isomorphism
(\ref{eqn:explicitsplitting}).
\end{lemma}

\begin{proof}
The independence of the
generator $\phi\in H^{1}(M;\ZZ)$ is clear: if we replace
$\phi$ by $-\phi$ we have to replace  
$s_p(1)$ by $s_p(-1)=-s_p(1)$;
% =z-p(z)$, where $\phi(z) =1$, by $-z-p(-z)$ since
% $-\phi(-z)=1$.
hence, $a=\alpha(s_{p}(1))$ has to be replaced by $a^{-1}$.
Moreover, we have to replace $\Delta_{k}^{\phi_{\sigma}}(t)$ by
$\Delta_{k}^{\phi_{\sigma}}(t^{-1})$ and 
the claim follows.

Suppose now that we have two projections 
$p_1,p_2\co H_1(M;\ZZ)\to\tors(H_1(M;\ZZ))$. Then
there is morphism
$\psi\co\ZZ\to \tors(H_1(M;\ZZ))$ as in (\ref{eqn:p2=p1-psi})
such that 
\begin{displaymath}
\psi\circ\phi=p_2-p_1=s_1\circ\phi-s_2\circ\phi,
\end{displaymath}
because $s_i\circ\phi=\operatorname{Id }-p_i$.
In particular 
\begin{displaymath}
\psi=s_1-s_2.
\end{displaymath}
Let $\sigma_i\co\p M\to\U$, $i=1,2$, be given by
$\sigma_i := \sigma\circ p_i$ and denote 
$a_i :=\alpha(s_{p_i}(1))$, so that $\alpha(\psi(1))=a_1 a_2^{-1}$.
By (\ref{eqn:deltaothersplitting}) we get:
\begin{displaymath}
\Delta^{\phi_{\sigma_2}}_k(t) = 
\Delta^{\phi_{\sigma_1}}_k(\alpha(\psi(1))\, t)=
\Delta^{\phi_{\sigma_1}}_k(a_1 a_2^{-1}\, t),
\end{displaymath}
Putting $t= a_2\, s$, we get
$
\Delta^{\phi_{\sigma_2}}_k(a_2\, s)= \Delta^{\phi_{\sigma_1}}_k(a_1 \, s)
$
Hence the order of vanishing of   $\Delta^{\phi_{\sigma_i}}_k$
at $a_i$ is independent of $i$.
\end{proof}

\begin{definition}
\label{definition:beta} Following~\ref{definition:zero}, we define
$\beta_{\alpha}\co\CC[t^{\pm 1}]\to\CC$
to be the evaluation map  at $\alpha(s_{p}(1))\in\mathbb C$, i.e.
$\beta_{\alpha}(\eta(t))=\eta(s_{p}(1))\in\CC$ $\forall 
\eta(t)\in\CC[t^{\pm 1}]$.
\end{definition}

The previous lemma says that the evaluation and the order of
vanishing of  $\Delta^{\phi_{\sigma}}_k$ at $\beta_{\alpha}$ is
independent of the splitting of the first homology group. In
addition, we have
\begin{equation}
\label{eqn:alphabetaphi} 
\alpha=\beta_{\alpha}\circ \phi_\sigma\,,
\end{equation}
 i.e.\ $\alpha(\gamma) =
\beta_{\alpha}(\sigma(\gamma)t^{\phi(\gamma)}) =
\sigma(\gamma)\alpha(s_{p}(1))^{\phi(\gamma)}$ $\forall \gamma\in \p
M$.

\begin{example}\label{ex:2iemepartie}
Let $M$ be the torus bundle given in Example~\ref{ex:1erepartie}.
For every $\lambda\in\CC^{*}$ there are representations
$\alpha_{i}\co\p M\to\CC^{*}$, $i=1,\ldots,4$, given by
$\alpha_{i}(\mu)=\lambda$ and $\alpha_{i}\big|_{\ZZ/2\oplus\ZZ/2}
 = \sigma_{i}$. Now, $\alpha_{i}$ is a root of the Alexander
invariant if and only if $\Delta^{\phi_{\sigma_i}}(\lambda)=0$. More precisely,
$\alpha_{1}$ is a root of the Alexander invariant if and only if
$\lambda= 3\pm \sqrt 8$ and,
 for $i=2,3,4$, $\alpha_{i}$ is a
root of the Alexander invariant if and only if $\lambda=1$. Note that in each
case the root is a simple root.
\end{example}

\begin{definition}
We define $\CC_{\alpha}$ to be the $\pi_1(M)$-module $\CC$ with the
action induced by $\alpha$,
i.e. $\gamma\cdot x= \alpha(\gamma) x$ $\forall x\in\CC$ and
$\forall\gamma\in\p M$.
\end{definition}

We explain the motivation of this definition.
Let $\rho_{\alpha}\co\p M\to\PSL(\CC)$ be the representation in
(\ref{eq:rho-alpha}). The Lie algebra $\Sl(\CC)$ turns into a
$\p M$-module via $\Ad\circ\rho_{\alpha}$ which will simply  be
denoted by $\Sl(\CC)_{\alpha}$. The  $\p M$-module
$\Sl(\CC)_{\alpha}$ decomposes as
$\Sl(\CC)_{\alpha}=\CC_+\oplus\CC_0\oplus\CC_-$, where
\begin{equation}\label{eqn:decomp}
\CC_{+}=\CC
\begin{pmatrix} 0 & 1\\ 0 & 0
\end{pmatrix},\quad
\CC_{0}=\CC
\begin{pmatrix} 1 & 0\\ 0 & -1
\end{pmatrix}
\text{ and }
\CC_{-}=\CC
\begin{pmatrix} 0 & 0\\ 1 & 0
\end{pmatrix}\,.
\end{equation}
Here $\CC_0$ is a trivial $\p M$-module and $\CC_{\pm}:=
\CC_{\alpha^{\pm}}$, where $\alpha^+=\alpha$ and
$\alpha^-$ is the morphism that maps every element $\gamma\in\pi_1(M)$
to $\alpha(\gamma^{-1})$.

\paragraph{Computing $\mathbf {H^{1}(\Gamma;\CC_{\alpha})}$}
Let $\Gamma = \pres{\xvec S n}{\xvec R m}$ be a finitely presented
group and let $\alpha\co\Gamma\to\CC^{*}$ be a representation.
In order to compute $H^{1}(\Gamma;\CC_{\alpha})$
we can use
the canonical 2-complex with one $0$-cell associated to the
presentation of $\Gamma$ (see Lemma~\ref{lem:EilMac}).
More precisely, we can identify
$Z^{1}(\Gamma;\CC_{\alpha})$ with the kernel of the linear map
$\CC^{n}\to\CC^{m}$ given by $\bfa \mapsto \bfA\, \bfa$ where
$\bfA=(a_{ji})\in M_{m,n}(\CC)$ is given by $a_{ji}=
\alpha(\partial R_{j}/ \partial S_{i})$.

For the remainder of this section we shall fix a projection 
$p\co H_{1}(M;\ZZ)\to\tors(H_{1}(M;\ZZ))$ and a generator
$\phi\in H^1(M;\ZZ)$.
Moreover, let $\sigma$ denote the
restriction of $\alpha$ to the torsion subgroup
$\tors(H_{1}(M;\ZZ))$.
Using (\ref{eqn:alphabetaphi})
 we obtain $a_{ji}= \beta_{\alpha} (J^{\phi_\sigma}_{ji})$ where
$J^{\phi_\sigma}=(J^{\phi_\sigma}_{ji})\in
M_{n-1,n}(\CC[t^{\pm 1}])$ is the Jacobian, i.e.\
$J^{\phi_\sigma}_{ji}:=\phi_{\sigma}\circ\pi(
\partial R_{j}/
\partial S_{i})$. We have $\dim B^{1}(\Gamma;\CC_{\alpha}) =1$ if
$\alpha$ is nontrivial, and $\dim B^{1}(\Gamma;\CC_{\alpha}) =0$ if
$\alpha$ is trivial. Hence for any nontrivial $\alpha\co\Gamma\to\CC^{*}$,
we have
\begin{equation}
\label{eq:dimH1}
\dim H^{1}(\Gamma;\CC_{\alpha}) = n - \Rk \bfA -1\,.
\end{equation}

\begin{lemma}\label{lem:dim} Let $\alpha\co\p M\to\CC^{*}$ be a
nontrivial homomorphism. Then $\dim H^{1}(\p M;\CC_{\alpha})=k$
if and only if $\alpha$
 is a zero of the $k$-th Alexander invariant and not a zero of the
$(k+1)$-th Alexander invariant.
\end{lemma}
\begin{proof} By (\ref{eq:dimH1}) we have $k=\dim
H^{1}(\p M;\CC_{\alpha}) = n - \Rk \bfA -1$, where
$\bfA=(\alpha(\partial R_{j}/
\partial S_{i}))$. Now $\alpha = \beta_{\alpha}\circ\phi_{\sigma}$ and hence
$\beta_{\alpha} (\Delta_{l}^{\phi_{\sigma}}(M)) = 0$ if $l<k$ and
$\beta_{\alpha} (\Delta_{l}^{\phi_{\sigma}}(M)) \neq 0$ if $l\geq k$.
\end{proof}
%\marnote{J'ai remplac\'e $\Delta_{l}^{\phi_{\sigma}}$ par
%$\Delta_{l}^{\phi_{\sigma}}(M)$.}

Recall that $\alpha^{\pm}\co\p M\to\CC^{*}$ denotes the homomorphisms
given by
\begin{equation}\label{eq:alpha/pm}
\alpha^{\pm}(\gamma) := (\alpha(\gamma))^{\pm 1}\qquad\forall\gamma\in\pi_1(M).
\end{equation}
It follows from the symmetry of the Alexander invariants that
$\alpha^{+}$ and $\alpha^{-}$ are zeros of the same order of the
$k$-th Alexander invariant. More precisely, we have:

\begin{prop}\label{prop:sameorder} Let
$\alpha\co\p M\to\CC^{\ast}$ be a nontrivial homomorphism. Then
$\alpha^{+}$ and $\alpha^{-}$ are zeros of the same order of the
$k$-th Alexander invariant. In particular, we have
$$
\dim H^{1}(\p M;\CC_{\alpha^{+}}) = \dim H^{1}(\p
M;\CC_{\alpha^{-}})\,.
$$ 
\end{prop}

\begin{proof}
Since $\alpha^{-}(s_{p}(1))$ is the inverse of $\alpha(s_{p}(1))$, it
suffices to check that $\Delta_k^{\phi_{\sigma^{-1}}}(t^{-1})=
\epsilon\,\Delta_k^{\phi_{\sigma}}(t)$ where $\epsilon$ is a unit
in $\CC[t^{\pm 1}]$. To verify this, notice that   the image of
$\sigma$ is contained in $\U$, so $\sigma^{-1}(
\gamma)=\overline{\sigma(\gamma)}$, $\forall\gamma\in\pi_1(M)$. Hence,
if $\Delta_k^{\phi_{\sigma}}(t)=\sum_i a_i t^i$, then
$\Delta_k^{\phi_{\sigma^{-1}}}(t)=\sum_i \overline{a_i} t^i$. By
Proposition~\ref{pro:symmetry},  $\sum_i a_i t^i$ differs from
$\sum_i \overline{a_i} t^{-i}$ by a unit, hence the proposition
follows.
\end{proof}

\paragraph{The space of abelian representations}
Let $\alpha\co\p M\to\CC^{*}$ be a representation and let
$\varphi\co\ZZ\to\CC^{*}$ be a homomorphism. Using multiplication, we
obtain a homomorphism $\alpha\varphi\co\p M\to \CC^{*}$ given by
$\alpha\varphi(\gamma)= \alpha(\gamma)\varphi(\phi(\gamma))$ for
$\gamma\in\p M$, where $\phi\in H^{1}(M;\ZZ)$ is a generator.
 There is  a one dimensional
irreducible algebraic set $V_{\alpha}\subset R(M)$ given by
$$V_{\alpha} := \{ \rho_{\alpha\varphi}\mid
\varphi\in\hom(\ZZ,\CC^{*})\}\subset R(M)\,.$$
Moreover,
the $\PSL(\CC)$ orbit of $\rho_{\alpha\varphi}$ is two dimensional if
$\alpha\varphi$
is nontrivial and hence $V_{\alpha}$ is contained in an at least three
dimensional component. We denote by
$S_{\alpha}(M)\subset R(M)$ the closure of the
$\PSL(\CC)$--orbit of $V_{\alpha}$.
Notice that $\rho_{\alpha}\in S_{\alpha}(M)$ and $\dim_{\CC}
S_{\alpha}(M)\geq 3$.

\begin{lemma} \label{lem:nondef}
Let $\alpha\co \p M \to \CC^{*}$ be a representation. If $\alpha$ is
not a zero of the Alexander invariant then there exists a neighborhood
of $\rho_{\alpha}$  in $R(M)$ consisting entirely of points of the component
$S_{\alpha}(M)$. Moreover, $\rho_{\alpha}\in S_{\alpha}(M)$ is a
smooth point and $S_{\alpha}(M)$ is the unique component through
$\rho_{\alpha}$ and $\dim S_{\alpha}(M) = 3$.
\end{lemma}
\begin{proof}
We have
$
3 \leq \dim S_{\alpha}(M)\leq \dim_{\rho_{\alpha}} R(M) \leq \dim
Z^{1}(\p M;\Sl(\CC)_{\alpha})
$.
Therefore Lemma~\ref{lem:locdim} implies the result if we can show
that 
\[ 
\dim Z^{1}(\p M;  \Sl(\CC)_{\alpha})=3\,.
\]

If $\alpha$ is trivial then
$Z^{1}(\p M;  \Sl(\CC)_{\alpha})=H^{1}(\p M;  \Sl(\CC)_{\alpha})\cong
H^{1}( M; \CC) \otimes \Sl(\CC)$ and
$\dim Z^{1}(\p M; \Sl(\CC)_{\alpha})=3$ follows.

If $\alpha$ is nontrivial, then the $\p M$-module
$\Sl(\CC)_{\rho_{\alpha}}$ splits as $ \Sl(\CC)_{\alpha} =
\CC_{+}\oplus\CC_{-}\oplus\CC_{0}$. Hence $H^{1}(\p M;
\Sl(\CC)_{\alpha})\cong H^{1}(\p M;\CC_{+})\oplus H^{1}(\p
M;\CC_{0})\oplus H^{1}(\p M;\CC_{-})$ and by Lemma~\ref{lem:dim}
and Proposition~\ref{prop:sameorder} we get 
\[ 
\dim H^{1}(\p M;\Sl(\CC)_{\alpha})= 1\,.
\] 
This implies that $\dim Z^{1}(\p
M;\Sl(\CC)_{\alpha})= 3$ and $\rho_{\alpha}\in
S_{\alpha}(M)\subset R(M)$ is a smooth point.
\end{proof}

\section{Cohomology of metabelian representations}
\label{sec:cohomologymetabelian}

The aim of the following three sections is to prove that, when
$\alpha$ is a simple zero of the Alexander invariant,
 certain reducible metabelian representations are smooth points of the
representation variety $R(M)$. First we construct these reducible
representations and then, before proving their smoothness in
Section~\ref{section:deformingmetabelian}, we  shall do some
cohomological  computations  in this section and the following
one.

Let  $\alpha\co\p M\to\CC^{*}$ be a homomorphism and let
$d\co\p M\to\CC$ be a map. The map 
$\rho^{d}_{\alpha}\co\p M\to\PSL(\CC)$ given by $$
\rho^{d}_{\alpha}(\gamma)=
\begin{pmatrix}
1 &  d(\gamma)\\
0 & 1
\end{pmatrix} \rho_{\alpha}(\gamma)=\pm
\begin{pmatrix}
\alpha^{\frac12}(\gamma) & \alpha^{-\frac12}(\gamma) d(\gamma)\\ 0 &
\alpha^{-\frac12}(\gamma)
\end{pmatrix}$$
is a homomorphism if and only if $d\in
Z^{1}(\p M;\CC_{\alpha})$. Moreover, $\rho^{d}_{\alpha}$ is non-abelian if and only if $d$ is not a coboundary.

\begin{cor}[Burde, de Rham] \label{cor:BurDeRha}
Let
$\alpha\co\p M\to\CC^{*}$ be a
representation and define $\rho_{\alpha}\co \p M\to\PSL(\CC)$ as in
(\ref{eq:rho-alpha}). Then there exists a reducible, non-abelian
representation $\rho\co\p M\to\PSL(\CC)$ such that
$\chi_{\rho}=\chi_{\rho_{\alpha}}$ if and only if
$\alpha$ is a zero of the Alexander invariant.
\end{cor}
\begin{proof}
By Lemma~\ref{lem:dim} we have that
$ \dim H^{1}(\p M;\CC_{\alpha})>0 $ if and only if
$\alpha$ is a zero of the  Alexander invariant of $M$.
\end{proof}

If $\alpha$ is a simple zero of
the Alexander invariant,
 then $\alpha^{\pm}$ defined by  (\ref{eq:alpha/pm})
 is a  zero of the first Alexander invariant,
but it is not a zero of the second. By Lemma~\ref{lem:dim} we have
 $H^1(\p M;\CC_{\pm})\cong\CC$.

Let  $d_{\pm}\in Z^1(\p M; \CC_{\pm})$ be a  cocycle which
represents a generator of the first cohomology group
$H^1(\p M; \CC_{\pm})$. We denote by
$\rho^{\pm}$  the metabelian representations into the
upper/lower triangular group given by Corollary~\ref{cor:BurDeRha}, i.e.\
$$
\rho^+(\gamma)=
\begin{pmatrix} 1 & d_+(\gamma)\\ 0 & 1
\end{pmatrix} \rho_\alpha(\gamma)
\qquad\textrm{ and }\qquad
\rho^-(\gamma)=
\begin{pmatrix} 1 & 0\\ d_-(\gamma) & 1
\end{pmatrix} \rho_\alpha(\gamma) . $$
If we replace  $d_{\pm}$ by $d'_{\pm}=c\,d_{\pm} + b_{\pm}$ where
$c\in\CC^{*}$ and $b_{\pm}\in B^{1}(\p M; \CC_{\pm})$ then
$\rho^{\pm}$ changes by conjugation by an upper/lower triangular
matrix.
Let $\mathfrak{b}_+ \subset \Sl(\CC)$ denote the Borel subalgebra of
upper triangular matrices. It is a $\p M$-module via
$\Ad\circ\rho^+$.
The short exact sequence
of $\p M$-modules $$
0\to\CC_+\to\mathfrak b_+\to\CC_0\to 0 $$
gives a long exact sequence in cohomology:
\begin{multline}\label{seq:long} 0\to
H^{0}(M;\CC_{0})\stackrel{\delta^{1}}{\longrightarrow}
H^{1}(M;\CC_{+})
\to H^{1}(M;\mathfrak b_+)\to\\ H^{1}(M;\CC_{0})
\stackrel{\delta^{2}}{\longrightarrow}H^{2}(M;\CC_{+})
\to H^{2}(M;\mathfrak b_+)\to 0
\end{multline}

\begin{lemma}
\label{lemma:b+acyclic}
 We have that $H^{1}(M;\mathfrak b_+)=0$ if and only if
$\delta^{2}\co H^{1}(M;\CC_{0})\to H^{2}(M;\CC_{+})$ is an
isomorphism.
\end{lemma}
\begin{proof} The Euler characteristic $\chi(M)$
vanishes.
Hence, $H^{1}(M;\mathfrak b_+)=0$ implies
$H^{2}(M;\mathfrak b_+)=0$ and  the sequence (\ref{seq:long}) gives
that $\delta^{2}\co H^{1}(M;\CC_{0})\to H^{2}(M;\CC_{+})$ is an
isomorphism.

Suppose that $\delta^{2}\co H^{1}(M;\CC_{0})\to H^{2}(M;\CC_{+})$ is
an isomorphism. Then the sequence (\ref{seq:long}) gives
$H^{2}(M;\mathfrak b_+)=0$ and the vanishing of the Euler
characteristic  implies $H^{1}(M;\mathfrak b_+)=0$.
\end{proof}

A cocycle $d_{0}\co\p M\to\CC_{0}$ is nothing but a homomorphism
$d_{0}\co\p M\to(\CC,+)$. A direct calculation gives
$$
\delta^{2}(d_{0})(\gamma_{1},\gamma_{2})= -2 d_{+}(\gamma_{1})
d_{0}(\gamma_{2}).
$$
The 2-cocycle $\delta^{2}(d_{0})$ is a cup product.
In our situation we
have the multiplication
$\CC_{0}\otimes\CC_{+}\to\CC_{+}$. This gives us a cup product
$$
\overset{.}{\cup}\co H^{1}(\p M;\CC_{0})\otimes
H^{1}(\p M;\CC_{+})\to H^{2}(\p M;\CC_{+})$$
and
\begin{equation}
\label{eqn:delta2iscup}
\delta^{2}(d_{0})= -2 (d_{+}\overset{.}{\cup}d_{0})
\end{equation}
 (see
Equation~(\ref{eq:cupprod})).
Hence we have that $\delta^{2}$ is an isomorphism if and only if,
for each nontrivial homomorphism $d_{0}\co\Gamma\to\CC$, the
cocycle $d_{+}\overset{.}{\cup}d_{0}$ represents a nontrivial
cohomology class.

The next lemma will be proved in Section~\ref{sec:cohomology}:
\begin{lemma} \label{lemma:nonzeroproduct}
Let $d_{0}\co\p M\to\CC_{0}$ be a nontrivial homomorphism and
let $d_{+}\co\p M\to\CC_{+}$ be a
cocycle representing a nontrivial cohomology class.
If $\alpha$ is
a simple zero of the Alexander invariant, then the 2-cocycle
$d_{+}\cupd d_{0}\in Z^{2}(\p M;\CC_{+})$ represents a nontrivial
cohomology class.
\end{lemma}

\begin{cor}\label{cor:dimh1}
Let $\alpha\co\p M\to\CC^{*}$ be a nontrivial
homomorphism. If $\alpha$ is a
simple zero of the Alexander invariant then
$H^1(\p M;\mathfrak{b}_+)=0$  and
the projection to the quotient
$\Sl(\CC)_{\rho^{+}}\to\Sl(\CC)_{\rho^{+}}/\mathfrak{b}_{+}\cong
\CC_-$ induces an isomorphism $$
H^1(\p M ;\Sl(\CC)_{\rho^+})\cong H^1(\p M;\CC_-)\cong \CC .$$
\end{cor}
\proof
Lemmas~\ref{lemma:nonzeroproduct} and
\ref{lemma:b+acyclic} and equation
(\ref{eqn:delta2iscup})
 imply that
$H^1(\p M;\mathfrak{b}_+)=0$. The isomorphism follows then from
 the long exact sequence in cohomology  corresponding to
$$
0\to \mathfrak b_+\to \Sl(\CC)_{\rho^{+}}\to\CC_-\to 0\,.\eqno{\qed} $$

\section{Fox calculus and 2-cocycles}
\label{sec:cohomology}

 The aim of this section is to prove Lemma~\ref{lemma:nonzeroproduct}.
 Let $\Gamma$ be a finitely presented group and let $A$ be a $\CC\Gamma$-module.
In the sequel we have to decide
when a given 2-\emph{cocycle} $c\co\Gamma\times\Gamma\to A$
is a coboundary.

A normalized 2-cochain is a map
$c\co\Gamma\times \Gamma\to A$ where the normalization condition is
$c(1,\gamma)= c(\gamma,1)=0$ for all $\gamma\in\Gamma$.
We shall extend $c$ linearly on the first component, i.e.\ for
$\eta = \sum_{\gamma\in\Gamma} n_{\gamma}\gamma\in\CC\Gamma$,
$n_{\gamma}\in\CC$, we define
$$
c (\eta,\gamma_{0}) := \sum_{\gamma\in\Gamma} n_{\gamma}\, c(\gamma,\gamma_{0})
\,.$$

\begin{prop}\label{pro:cobo}
        Let $\Gamma = \pres{\xvec{S}{n}}{\xvec{R}{m}}$ be a finitely
        presented group and let $c\co\Gamma\times \Gamma\to A$ be
        normalized 2-cocycle.

        Then $c\co\Gamma\times\Gamma\to A$ is a
        coboundary if and only if there exists $a_{i}\in A$,
        $i=1,\ldots,n$, such that for all $j=1,\ldots,m$ the equation
\begin{equation}
\label{eq:cobound}
                        \sum_{i=1}^n \pi\big(
                \frac{
\partial R_{j}}{
\partial S_i} \big) \circ a_i
                             + \sum_{i=1}^n
                 c (\pi(\frac{
\partial R_{j}}{
\partial S_i}),\pi(S_i)) =0
\end{equation}
         holds.
\end{prop}

Here $\pi\co\langle S_1,\ldots,S_n\rangle\to \Gamma$
denotes the
natural projection from the free group to $\Gamma$.
\begin{proof}
We start by recalling some well known constructions used
in the proof (cf.~\cite{Bro82}).
Let $X$ be the canonical 2-complex with one $0$-cell
associated
to the presentation of $\Gamma=\pres{\xvec{S}{n}}{\xvec{R}{m}}$
i.e.\ $X = X^{0}\cup X^{1}\cup X^{2}$ where $X^{0}=\{e^{0}\}$,
$X^{1}= \{ \xvec {e^{1}} n \}$ and $X^{2}= \{ \xvec {e^{2}} m \}$.
The universal covering $p\co\tilde X\to X$ gives us a free chain
complex $C_{k} :=C_{k}(\tilde X)$ of $\Gamma$-modules. A basis
for $C_{k}$ is given by choosing exactly one cell
$\tilde e^{k}_{j}\in p^{-1}(e^{k}_{j})$.
With respect to this basis, $\partial_{2}$ is given by the
Fox calculus, i.e.\ 
\begin{equation}\label{eqn:boundaries}
\partial_{2}(\tilde e^{2}_{j}) =
\sum_{i=1}^{n} \frac{\partial R_{j}}{\partial S_{i}} \tilde e^1_i
\quad\text{ and}\quad
\partial_{1} (\tilde e^{1}_{i})= (S_{i}-1) \tilde e^{0}\,.
\end{equation}
Notice that $\partial_{1}\circ\partial_{2}=0$ corresponds to the
fundamental formula of the Fox calculus (see \cite{BuZi03}).

The  \emph{normalized} bar resolution for $\Gamma$ is denoted by
$B_{*}:=B_{*}(\Gamma)$. More precisely, let
$B_{n} := B_{n}(\Gamma)$ be the free $\Gamma$-module with
generators $[x_{1}|\ldots | x_{n}]$, where
$x_{i}\in\Gamma\setminus\{1\}$. In order to give meaning to every
symbol $[x_{1}|\ldots | x_{n}]$ set
$[x_{1}|\ldots | x_{n}] = 0$ if $x_{i} =1$ for any $i$.
This is called the \emph{normalization} condition.
Note that $B_{0}\cong\ZZ\Gamma$
is the free $\Gamma$-module on one generator  and  the
augmentation $\varepsilon\co B_{0}\to\ZZ$ maps $[\, ]$ to $1$.
In low dimensions the boundary operators  are given by
\[
\partial [x|y]= x[y] -[xy] + [x],\quad
\partial[x]= (x-1)[\,] \,.
\]
Moreover, homomorphisms $s_{-1}\co\ZZ\to B_{0}$, $s_{n}\co B_{n}\to
B_{n+1}$ of abelian groups are defined by $$
s_{-1}(1) = [\;] \text{ and }
s_{n}(x[x_{1}|\ldots | x_{n}])=[x|x_{1}|\ldots | x_{n}]\;.$$
It turns out that $(s_{n})$ is a
contracting homotopy for the underlying augmented chain
complex $B_{*}\stackrel{\varepsilon}{\to}\ZZ$ of
abelian groups.
$$
\begin{CD}
 0 @>>> C_{2}(\tilde X) @>\partial_{2}>> C_{1}(\tilde X)
@>\partial_{0}>> C_{0}(\tilde X) @>\epsilon>> \ZZ\\
@. @VVf_{2}V @VVf_{1}V @VVf_{0}V @VV=V\\
\cdots @>>>
B_{2}(\Gamma) @>\partial>\stackrel{\longleftarrow}{s_{1}}>
B_{1}(\Gamma) @>\partial>\stackrel{\longleftarrow}{s_{0}}>
B_{0}(\Gamma) @>\epsilon>\stackrel{\longleftarrow}{s_{-1}}> \ZZ
\end{CD}
$$
By using the fact that
$C_{*}$ is a  free $\Gamma$-complex and that $B_{*}$ is contractile
we obtain a chain map $f _{*}\co C_{*}\to B_{*}$ which
is augmentation preserving i.e.\ $\varepsilon \circ f_{0} =
\varepsilon$. Moreover, $f$ is unique up to chain homotopy (see
\cite[I.7.4]{Bro82}). The contracting homotopy $s_{n}\co B_{n}\to
B_{n+1}$ and the basis of $C_{*}$ determine $f_{*}$
inductively (see \cite[p.24]{Bro82} for details).
Hence the maps $f_{i}\co C_{i}\to B_{i}$, $i=0,1,2$, are given by
\begin{equation}\label{eqn:f}
f_{0}(\tilde e^{0})= [\;], \quad
f_{1}(\tilde e^{1}_{i})=[S_{i}] \text{ and }
 f_{2}(\tilde e^{2}_{j})=\sum_{i=1}^{n}[\frac{
\partial R_{j}}{
\partial S_{i}}\big| S_{i}]\,.
\end{equation}
Here we have used the following convention:  if
$\eta=\sum_{\gamma\in\Gamma}
n_{\gamma}\gamma\in\ZZ\Gamma$, then
$[\eta,\gamma_{0}]:= \sum_{\gamma\in\Gamma}
n_{\gamma}[\gamma,\gamma_{0}]\in B_{2}$.

It follows that the induced cochain map
$f_{2}^{*}\co \hom_{\Gamma}(B_{2};A)\to\hom_{\Gamma}(C_{2};A)$
is given by
\begin{equation}\label{eqn:f2}
f_{2}^{*}(c) (\tilde e^{2}_{j}) =
\sum_{i=1}^{n} c\big(\frac{
\partial R_{j}}{
\partial S_{i}}, S_{i}\big), \text{ where $c\in
\hom_{\Gamma}(B_{2};A)$.}
\end{equation}
By
Lemma~\ref{lem:EilMac}, the map $f_{2}^{*}$
induces an injection $f^{*}\co H^{2}(\Gamma ; A)\to H^{2}(X; A) $.
 Hence, $f_{2}^{*}(c)$ is a coboundary if and only if
there exists a cochain 
\begin{equation}\label{eqn:lastf2}
b\in\hom_{\Gamma}(B_{1};A) \text{ such that }
f_{2}^{*}(c) (\tilde e^{2}_{j}) + b(
f_{1}(\partial_{2} \tilde e^{2}_{j} )) = 0
\end{equation}
for all $j=1,\ldots,m$.
The proposition follows from Equation~(\ref{eqn:lastf2}) by using
 (\ref{eqn:boundaries}), (\ref{eqn:f}) and (\ref{eqn:f2}).
\end{proof}

Let $\alpha\co\p M\to\CC^{*}$ be a nontrivial representation.
Note that $\alpha$ induces a homomorphism $\alpha\co
\tors(H_{1}(M;\ZZ))\oplus\ZZ \to \CC^{*}$, and
let $\sigma$ denote the
restriction of $\alpha$ to the torsion subgroup
$\tors(H_{1}(M;\ZZ))$.
As before we
denote by $\CC_{+}$ the $\p M$-module $\CC_\alpha$, i.e.\ $\gamma\cdot
z=\alpha(\gamma) z$. For the remainder of the section, fix a
projection $p\co H_{1}(M;\ZZ)\to\tors(H_{1}(M;\ZZ))$ and a generator
$\phi\in H^1(M;\ZZ)$.

Let $h\co\p M\to(\CC,+)$ be a
homomorphism and let $d_{+}\co\p M\to\CC_{+}$ be a cocycle.
Since $\phi$ is a generator of the first cohomology group, there exists
$a\in\CC$ such that $h = a\, \phi$.

By Proposition~\ref{pro:cobo} we have that $h\overset{.}{\cup}d_{+}$
is a coboundary if and only if there exist
$x_{1},\ldots,x_{n}\in\CC_{+}$ such that
\begin{equation}
\label{eq:hcupd}
\sum_{i=1}^{n}
\alpha (\frac{
\partial R_{j}}{
\partial S_{i}}) \,x_{i} +  h\overset{.}{\cup}d_{+}
(\frac{
\partial R_{j}}{
\partial S_{i}},S_{i}) = 0
\end{equation}
for all $j=1,\ldots,n-1$.

We have
$ h\overset{.}{\cup}d_{+} (\gamma_{1},\gamma_{2}) =
h(\gamma_{1})\alpha(\gamma_{1})d_{+}(\gamma_{2})$ and hence for
$\eta=\sum c_{\gamma}\gamma\in\CC\p M$ we get from (\ref{eqn:alphabetaphi}):
\begin{align} h\overset{.}{\cup}d_{+}(\eta,\gamma_{0}) &=
\sum c_{\gamma}\,(h\overset{.}{\cup}d_{+}) (\gamma,\gamma_{0})\notag \\
&= \label{equ:cup}
\sum c_{\gamma}\, h(\gamma)\alpha(\gamma)d_{+}(\gamma_{0})\notag \\
&= a \sum
c_{\gamma}\,\phi(\gamma)\sigma(\gamma)\alpha(s_{p}(1))^{\phi(\gamma)}d_{+}(\gamma_{0}
)\,.
\end{align}
%\marnote{Nouvelle ligne\ldots}

Let $D\co\CC[t^{\pm 1}]\to\CC[t^{\pm 1}]$ be the following differential
operator:
\[
D(\sum_{i\in\ZZ} c_{i}t^i) = \sum_{i\in\ZZ} c_{i} i t^i\,.
\]
The operator $D$ satisfies  the  following rules:
\begin{xalignat*}{2}
 D(c_{1}\eta_{1}+c_{2}\eta_{2}) &=
 c_{1}D(\eta_{1}) +  c_{2}D(\eta_{2})
 && \text{ for $c_{i}\in\CC$ and $\eta_{i}\in\CC[t^{\pm 1}]$,}\\
 D(\eta_{1}\eta_{2}) &= D(\eta_{1})\eta_{2} +
\eta_{1} D(\eta_{2})
 && \text{ for $\eta_{i}\in\CC[t^{\pm 1}]$.}
\end{xalignat*}
It follows from these rules that
$D(c)= 0$ for $c\in\CC$.

For given $z\in\CC^{*}$  we define $\ord_{z}\co\CC[t^{\pm 1}]\to
\NN\cup\{\infty\}$ to be the order of $\eta(t)\in \CC[t^{\pm 1}]$ at
$z$, i.e.\ $\ord_{z}(\eta) = \infty \Leftrightarrow \eta(t) \equiv 0$
and
\[
\ord_{z}(\eta) = k\in\NN \Leftrightarrow \exists \eta'\in \CC[t^{\pm 1}]
:\, \eta'(z)\neq 0  \text{ and }\eta(t) = (t-z)^{k} \eta'(t)\,.
\]
It is easy
to see that if $\eta(z)=0$, then
$\ord_{ z}(\eta)=\ord_{ z}(D(\eta))+1$.

For a fixed $z\in\CC^{*}$ the evaluation map $\CC[t^{\pm 1}]\to\CC$
which maps
$\eta(t)$ to $\eta(z)$ turns $\CC$ into a $\CC[t^{\pm 1}]$-module
which will be denoted by $\CC_{z}$.
The kernel of the evaluation map
$\CC[t^{\pm 1}]\to\CC$ is
exactly the maximal ideal generated by $(t-z)$.

We choose a splitting  of $H_1(M;\ZZ)$ as in
(\ref{eqn:explicitsplitting}) and we write
$\alpha=\beta_{\alpha}\circ\phi_{\sigma}$ as in
(\ref{eqn:alphabetaphi}). Recall that $\beta_{\alpha}$ is nothing but
the evaluation map at
$z :=\alpha(s_{p}(1))$. Hence (\ref{equ:cup}) gives:
\begin{equation}
\label{eq:der}
  h\overset{.}{\cup}d_{+} (\eta,\gamma_{0}) =
            a \,\beta_{\alpha}(D(\phi_{\sigma}(\eta))) d_{+}(\gamma_{0}),
\end{equation}
where $\eta\in\CC\p M$, $\phi_{\sigma}(\eta)\in\CC[t^{\pm 1}]$ and
$\gamma_{0}\in\p M$.

A $\CC[t^{\pm 1}]$-module homomorphism
$f\co(\CC[t^{\pm 1}])^{n}\to(\CC[t^{\pm 1}])^{m}$ induces a
$\CC[t^{\pm 1}]$-morphism
$f^z\co\CC^{n}_z\to\CC^{m}_z$.
This follows simply from $f(\Ker {(\beta_{\alpha})^{n}})\subset
\Ker {(\beta_{\alpha})^{m}}$.
\[
\begin{CD} (\CC[t^{\pm 1}])^{n} @>f>> (\CC[t^{\pm 1}])^{m}\\
@VV\beta_{\alpha}^{n}V  @VV\beta_{\alpha}^{m}V\\
\CC^{n}_z @>f^z>> \CC^{m}_z\,.
\end{CD}
\]
It is easy to see that
$D(f)\co(\CC[t^{\pm 1}])^{n}\to(\CC[t^{\pm 1}])^{m}$ given by
$D(f) := D^{m}\circ f - f \circ D^{n}$ is a
$\CC[t^{\pm 1}]$-module morphism.  If $\bfA$ is the matrix of $f$ with
respect to the canonical basis, then $\beta_{\alpha}(\bfA)$ is the
matrix of $f^z$ with
respect to the canonical basis and $D\bfA$ is the matrix of
$D(f)$ with
respect to the canonical basis. Here, $\beta_{\alpha}$ and $D$ applied
to a matrix means simply applying it to each entry.

\begin{proof}[Proof of Lemma~\ref{lemma:nonzeroproduct}]
Recall that we made the assumption  that $\alpha$
is a simple zero of the Alexander invariant of $M$.

Let $h\in Z^{1}(\p M;\CC_{0})$ and $d_{+}\in Z^{1}(\p M;\CC_{+})$ be
cocycles representing nontrivial cohomology classes.
It follows from Equation~(\ref{eq:cupsym}) that
\[
h\overset{.}{\cup}d_{+} + d_{+}\overset{.}{\cup}h\sim 0,
 \]
and hence $h\overset{.}{\cup}d_{+} $ is a coboundary if and only if
$d_{+}\overset{.}{\cup}h$ is a coboundary.

Let $a\in\CC^{*}$ such that $h=a\phi$ and set $z:= \alpha(s_{p}(1))$.
Then $\beta_{\alpha}(\eta(t))$ is simply $\eta(z)$.
By writing
equations~(\ref{eq:hcupd}) in matrix form, we
obtain from (\ref{eq:der}) that $h\overset.\cup d_{+}$ is a coboundary
if and only if the system
\begin{equation}
\label{eq:hcupd2}
J^{\phi_{\sigma}}(z)\, \bfx + a\, (DJ^{\phi_{\sigma}})(z)\,
\begin{pmatrix}
d_{+}(S_{1})\\ \vdots\\  d_{+}(S_{n})
\end{pmatrix}= 0
\end{equation}
has a solution $\bfx\in\CC^{n}$. Here, for each matrix $A\in
M_{m,n}(\CC[t^{\pm 1}])$ we denote by $A(z)\in M_{m,n}(\CC)$
the matrix obtained by applying the evaluation map to its entries.

From the canonical 2-complex associated to the presentation
 we
obtain the following resolutions:
$$
\begin{CD} 0 @<<<  (\CC[t^{\pm 1}])^{n-1} @<d_{2}<< (\CC[t^{\pm 1}])^{n}
@<d_{1}<<
\CC[t^{\pm 1}]
 @<<<  0\,\mbox{\,}\\
 @. @VV{\beta_{\alpha}}V   @VV{\beta_{\alpha}}V  @VV{\beta_{\alpha}}V @.\\ 0 @<<
<
\CC^{n-1}_z
@<{d_{2}^z}<<  \CC^{n}_z
@<{d_{1}^z}<< \CC_z @<<< 0\,.
\end{CD}
$$
The matrix of $d_{2}$ (respectively $d_{2}^z$)
with respect to the canonical basis is $J^{\phi_{\sigma}}$ (respectively
$J^{\phi_{\sigma}}(z)$) and
the matrix of $D(d_{2})$ with respect to the canonical basis is
$DJ^{\phi_{\sigma}}$.

It follows that $h\overset.\cup d_{+}$ is a coboundary if and only if
\[
(DJ^{\phi_{\sigma}})(z)
\begin{pmatrix} d_{+}(S_{1})\\ \vdots\\  d_{+}(S_{n})
\end{pmatrix} \in \image (d_{2}^z)\,.
\]
Note that $J^{\phi_{\sigma}}$ is a presentation matrix of
$H_{1}^{\phi_{\sigma}}(M,x_{0})\cong
H_{1}^{\phi_{\sigma}}(M)\oplus\CC[t^{\pm 1}]$.
The assumption that $\alpha$ is a simple zero of the Alexander
invariant implies that $H_{1}^{\phi_{\sigma}}(M)$ is torsion.
Hence,
there exist a basis
$\mathcal B = (b_{0},\ldots,b_{n-1})$ of $(\CC[t^{\pm 1}])^{n}$ and a
basis
$\mathcal C = (c_{1},\ldots,c_{n-1})$ of $(\CC[t^{\pm 1}])^{n-1}$ such
that $d_2(b_{0}) = 0$ and $d_2(b_{i}) = r_{i}(t) c_{i}$, $1\leq i\leq
n-1$, where $r_{i}(t)\in\CC[t^{\pm 1}]$ are nonzero and $r_{i+1}(t)\mid
r_{i}(t)$. Moreover, we have that $\ord_z(\Delta^{\phi_{\sigma}})=1$
and hence $r_{1}(z)=0$ and $(D r_{1})(z)\neq
0$. In particular, $r_{j} (z)\neq 0$ for $j\geq 2$.
%\marnote{J'ai remplac\'e $\ord_z(\Delta_{1}^{\phi_{\sigma}})=1$ par
%$\ord_z(\Delta^{\phi_{\sigma}})=1$.}

Now, $d_{2}\circ d_{1} =0$ gives that $\image d_{1}\subset
\CC[t^{\pm 1}]\cdot b_{0}$. Therefore, there exists a
$r(t)\in\CC[t^{\pm 1}]$ such that $d_{1}(1)= r(t)  b_{0}$. Since
$H^{0}(M,\CC_{+})=0$ we obtain $r(z)\neq 0$.

We define a basis $(b_{0}^z, \ldots,b_{n-1}^z)$ of
$\CC^{n}_z$ by  $b_{0}^z = r(z)  \beta_{\alpha}(b_{0}) = r(z) b_{0}(z)$
and
$b_{i}^z = \beta_{\alpha} (b_{i}) = b_{i}(z)$, $1\leq i\leq n-1$.
Analogously, a basis $(c_{1}^z, \ldots,c_{n-1}^z)$ of $\CC^{n-1}$ is
given by
$c_{i}^z = \beta_{\alpha}(c_{i}) = c_{i}(z)$, $1\leq i\leq n-1$.

We have
\[
\image d_{2}^z =
\Span(c_{2}^z, \ldots,c_{n-1}^z)\text{ and }
\Ker d_{2}^z =
\Span(b_{0}^z,b_{1}^z)\,.
\]
Note that $\Ker d_{2}^z = \Span(b_{0}^z,b_{1}^z)$ can be identified with
$Z^{1}(M;\CC_+)$ and that the coboundaries correspond to the
multiples of $b_{0}^z$.

Now a direct calculation gives
$$
(D d_{2})(z)\, (b_{0}^z) =
\beta_{\alpha}( D(d_{2})\, (r(t) b_{0})) \in \image
d_{2}^z,  \text{ (using $d_{2}(b_{0})= 0$),}$$
$$(Dd_{2})(z)\,(b_{1}^z) =
 \beta_{\alpha}( D(d_{2}) (b_{1}) ) \in (D r_{1})(z)\, c_{1}^z + \image
d_{2}^z,
 \text{ (using $r_{1}(z)=0$).}\leqno{\hbox{and}}$$
Moreover, $(D r_{1})(z)\neq 0$ and each element of
$\Ker d_{2}^z$ representing a nonzero cohomology class
does not map into $\image d^{z}_{2}$ under $(D d_{2})(z)$.
Hence,
for each cocycle $d_{+}\co\p M\to\CC_{+}$ which
represents a generator of $H^{1}(M;\CC_{+})$ and each nontrivial
homomorphism $h\co \p M\to(\CC,+)$ the system (\ref{eq:hcupd2})
has no solution, i.e.\ $h\overset{.}{\cup} d_{+}$ is not a coboundary.
\end{proof}

\section{Deforming metabelian representations}
\label{section:deformingmetabelian} We  suppose in the sequel
that $\alpha$ is a simple zero of the Alexander invariant of $M$.
Let $\rho^+\in R(M)$ denote the metabelian representation defined
in Section~\ref{sec:cohomologymetabelian}.  In this section we use
the results of the previous two sections in order to show that
$\rho^+$ is a smooth point of $R(M)$ with local dimension 4.

Let $i\co
\partial M\to M$ be the inclusion.
\begin{lemma}\label{lemma:notrivialonboundary}
The representation
$\rho^+\circ i_{\#}\co \p {\partial M}\to\PSL(\CC)$ is nontrivial.
\end{lemma}

\begin{proof} By   Corollary~\ref{cor:dimh1}, $H^1(
M;\Sl(\CC)_{\rho^+})\cong H^1(\p M;\Sl(\CC)_{\rho^+})\cong
\CC$ and by duality, $H^2(M;\partial
M;{\Sl(\CC)}_{\rho^+})\cong H^1( M;{\Sl(\CC)}_{\rho^+})\cong \CC$.
Thus, by the exact sequence of the pair $$
   H^1(M;{\Sl(\CC)}_{\rho^+})\to  H^1(
\partial M;{\Sl(\CC)}_{\rho^+}) \to H^2(M;
\partial M;{\Sl(\CC)}_{\rho^+}), $$
we obtain $\dim H^1(
\partial M;{\Sl(\CC)}_{\rho^+})\leq 2$.

We prove the lemma by contradiction:
if $\rho^+$ restricted to $i_{\#}( \p {\partial M})$ was trivial,
then ${\Sl(\CC)}_{\rho^+}$ would be a trivial $\p{\partial
M}$-module, and therefore $$ H^1(
\partial M;{\Sl(\CC)}_{\rho^+})\cong H^1(
\partial M;\CC)\otimes_\CC{\Sl(\CC)}\cong
{\Sl(\CC)}\oplus{\Sl(\CC)} $$ would have dimension six,
contradicting the previous upper bound for the dimension.
\end{proof}

\begin{definition}
  \label{def:Klein}
  A non-cyclic  abelian subgroup of $\PSL(\CC)$ with four elements is
  called \emph{ Klein's 4-group}. Such a group is realized by
  rotations about three orthogonal geodesics  and it is conjugate
  to the one generated by
  $\pm\left(
\begin{smallmatrix} 0 & 1\\ -1 & 0
\end{smallmatrix} \right)$ and
  $\pm \left(
\begin{smallmatrix} i & 0\\ 0 & -i
\end{smallmatrix} \right)$.
\end{definition}

\begin{remark}\label{rem:nokleinonboundary}
 The image $\rho^+ (i_{\#}(\p {\partial M}))$
cannot be the Klein group, because the image of $\rho^+$ is reducible
(i.e. the action on $ P^1(\CC)=\CC\cup\{\infty\}$
 has
a fixed point, $\infty$), but the Klein group has no fixed point in
$ P^1(\CC)$.
\end{remark}
%\marnote{J'ai remplace $\mathbf P^{1}(\CC)$ par $P^{1}(\CC)$ pour
%\^{e}tre compatible avec l'introduction.}

Recall that by Weil's construction $Z^1(\p M;{\Sl(\CC)}_{\rho^+})$
contains the Zariski tangent space of $R(M)$ at $\rho^+$ (cf.\
Section~\ref{sec:notation}). To prove the smoothness, we show that
all cocycles in $Z^1(\p M;{\Sl(\CC)}_{\rho^+})\cong \CC$ are
integrable. To do this, we prove that all obstructions vanish, by
using the fact that the obstructions vanish on the boundary.

\begin{lemma}\label{lemma:torussmooth} The variety $R(\ZZ\oplus\ZZ)$
has exactly two irreducible components. One is four dimensional and
smooth except at the trivial representation. The other component is
three dimensional and smooth; it is  exactly the orbit of a
representation onto the Klein group.
\end{lemma}

\begin{proof}
Let $\ZZ\oplus\ZZ=\pres{x,y}{[x,y]=1}$ and let
$\rho\co\ZZ\oplus\ZZ\to\PSL(\CC)$ be a representation given by
$\rho(x)=
\pm A_{x}$ and $\rho(y)= \pm A_{y}$. Then $\tr [A_{x},A_{y}]\in\{\pm
2\}$ and $\rho$ lifts to $\SL(\CC)$ if and only if $\tr
[A_{x},A_{y}]=2$.
Moreover, $\rho$ is a representation onto a Klein group if and only if
$\tr [A_{x},A_{y}]=-2$.

Thus $R(\ZZ\oplus\ZZ)$ has two components, one of dimension four and
one of dimension three, which is the orbit of a representation onto
the Klein group.

Given a representation  $\varrho\in R(\ZZ\oplus\ZZ)$ which is nontrivial and different from the Klein group, then
$H^0(\ZZ\oplus\ZZ;{\Sl(\CC)}_{\varrho}) \cong
{\Sl(\CC)}^{\varrho(\ZZ\oplus\ZZ)}\cong\CC$. Thus, by duality and
Euler characteristic, $H^1(\ZZ\oplus\ZZ;{\Sl(\CC)}_{\varrho})\cong
\CC^2$ and $Z^1(\ZZ\oplus\ZZ;{\Sl(\CC)}_{\varrho})\cong \CC^4$. This
computation shows that the dimension of the Zariski tangent space at
this representation is at most four. Since the representation lies in
a four dimensional component, it is a smooth point of
$R(\ZZ\oplus\ZZ)$. Note that it follows from the proof of 
Lemma~\ref{lemma:notrivialonboundary} that the trivial representation is a singular point of $R(\ZZ\oplus\ZZ)$.

If $\varrho$ is a representation onto the Klein group then
$H^0(\ZZ\oplus\ZZ;{\Sl(\CC)}_{\varrho})=0$ and hence
$H^1(\ZZ\oplus\ZZ;{\Sl(\CC)}_{\varrho})=0$ by the same Euler
characteristic argument. Hence
$\dim Z^1(\ZZ\oplus\ZZ;{\Sl(\CC)}_{\varrho})=3$. Since the orbit of
the representation is three dimensional and closed, the lemma is
proved.
\end{proof}

Given a cocycle $Z^1(\p M;{\Sl(\CC)}_{\rho^+})$ the first
obstruction to integration is the cup
product with itself. In general when the $n$-th obstruction vanishes,
the obstruction of order $n+1$ is defined, it lives in
$H^2(\p M;{\Sl(\CC)}_{\rho^+})$.

Let $\Gamma$ be a  finitely presented group and let
$\rho\co \Gamma \to \PSL(\CC)$ be a representation.
A \emph{formal deformation} of $\rho$ is a homomorphism
$\rho_{\infty}\co \Gamma\to\PSL(\CC [[t]])$
$$  \rho_{ \infty} ( \gamma) = \pm \exp ( \sum^{ \infty}_{i=1} t^i u_i(
        \gamma) ) \rho( \gamma)$$
where $u_i\co \Gamma\to\Sl(\CC)$ are
elements of $C^1( \Gamma,\Sl(\CC)_{\rho})$ such that
$ p_{0}\circ\rho_{ \infty}  = \rho$.
Here  $p_{0}\co\PSL(\CC [[t]])\to\PSL(\CC)$ is the evaluation
homomorphism
at $t=0 $ and  $ \CC [[t]] $ denotes  the ring of formal power
series. We shall say that  $\rho_{\infty}$ is a \emph{formal
deformation up to order $k$} of $\rho$ if $\rho_{\infty}$ is  a
homomorphism modulo $t^{k+1}$.

An easy calculation gives
that $\rho_{\infty}$ is a homomorphism up to first order if and only if
$u_1 \in Z^1(\Gamma,\Sl(\CC)_{\rho})$ is a cocycle.
We call a cocycle
$u_1 \in Z^1(
\Gamma,\Sl(\CC)_{\rho})$ \emph{formally integrable} if there is a
formal deformation of $\rho$ with leading term $u_1$.

Let $u_{1},\ldots,u_{k}\in C^1( \Gamma,\Sl(\CC)_{\rho}) $ such that
$$\rho_{k} (\gamma) = \exp ( \sum^{ k}_{i=1} t^i u_i(
        \gamma) ) \rho( \gamma)$$
 is a homomorphism into $\PSL(\CC [[t]])$ modulo $t^{k+1}$. 
Then there exists an
obstruction class $\zeta_{k+1} := \zeta_{k+1}^{(u_1,\ldots,u_{k})}
                                \in
H^2(\Gamma;\Sl(\CC)_{\rho})$
    with the following properties (see \cite[Sec.~3]{HPS01def}):
    \begin{description}
        \item[(i)]  There is a cochain $u_{k+1}\co\Gamma\to\Sl(\CC)$
such that
        \[
        \rho_{k+1} (\gamma) = \exp ( \sum^{ k+1}_{i=1} t^i u_i(
        \gamma) ) \rho( \gamma)
        \]
         is a homomorphism modulo $t^{k+2}$
if and only if $\zeta_{k+1}=0$.

        \item[(ii)]  The obstruction $\zeta_{k+1}$ is natural, i.e.\
        if $f\co\Gamma'\to\Gamma$ is a homomorphism   then
        $f^{*}\rho_{k} :=  \rho_{k} \circ f $
        is also a homomorphism modulo $t^{k+1}$ and
        $f^{*} (\zeta_{k+1}^{(u_1,\ldots,u_{k})}) =
    \zeta_{k+1}^{(f^{*}u_1,\ldots,f^{*}u_{k})}$.
        \end{description}

\begin{lemma} \label{lemma:obstructionsvanish}
    Let
    $\rho\co\p M\to\PSL(\CC)$ be a reducible, nonabelian 
    representation such that
    $\dim Z^{1}(\p M;\Sl(\CC)_{\rho}) =4$.

    If $\rho\circ i_{\#}\co \p{\partial M}\to\PSL(\CC)$ is
    neither trivial nor a representation onto a Klein group, then
    every cocycle in $Z^{1}(\p M;\Sl(\CC)_{\rho})$ is integrable.
\end{lemma}

\begin{proof}
We first show that
$i^{*}\co H^2(\p M;\Sl(\CC)_{\rho})\to
H^2(\p{\partial M}; \Sl(\CC)_{\rho})$ is injective.
To that extent
we use the following commutative diagram: $$
\begin{CD}
    H^2(M;{\Sl(\CC)}_{\rho}) @>\cong >> H^2(
\partial M;{\Sl(\CC)}_{\rho}) \\
    @AAA         @AA\cong A \\
    H^2(\p M;{\Sl(\CC)}_{\rho}) @>>>
    H^2(\p{\partial M};{\Sl(\CC)}_{\rho})
\end{CD} $$
The horizontal isomorphism on the top of the diagram comes from the
exact sequence of the pair $(M,\partial M)$ and the dimension
computation in the proof of Lemma~\ref{lemma:notrivialonboundary}. The
vertical isomorphism on the right is a consequence of asphericity of
$\partial M$. In addition, the vertical map on the left
 is an injection (see Lemma~\ref{lem:EilMac}).

    We shall now prove  that every element of
    $Z^1(\p M;\Sl(\CC)_{\rho})$ is
    integrable.  Let $u_1,\ldots,u_k\co\p M\to\Sl(\CC)$
    be given such that
    $\rho_{k} (\gamma) = \exp ( \sum^{ k}_{i=1} t^i u_i(
        \gamma) ) \rho( \gamma)$
 is a homomorphism modulo $t^{k+1}$.
Then the restriction
$\rho_{k}\circ i_{\#}\co \p{\partial M}\to\SL(\CC[[t]])$
is also a formal deformation of order $k$. On the other hand, it
follows from Lemma~\ref{lemma:torussmooth} that
 the restriction $ \rho_{k}\circ i_{\#}$ is a smooth point of the
representation variety $R(\partial M)$.
Hence,  the formal implicit function theorem gives that
$i^{*}\rho_{k}$ extends to a formal deformation of order $k+1$ (see
\cite[Lemma~3.7]{HPS01def}). Therefore, we have that
\[
0=\zeta_{k+1}^{(i^{*}u_1,\ldots,i^{*}u_{k})} = i^{*}
\zeta_{k+1}^{(u_1,\ldots,u_{k})}\,.
\]
Now, $i^{*}$ is injective and the obstruction vanishes.
\end{proof}

\begin{prop}
\label{prop:rho+smooth}
 The representation $\rho^+$ is a smooth point of $R(M)$ with local
dimension four.
\end{prop}

\begin{proof}
It follows from Corollary~\ref{cor:dimh1}
that $\dim H^1(\p M;{\Sl(\CC)}_{\rho^+})=1$.
Thus $\dim Z^1(\p M;{\Sl(\CC)}_{\rho^+})=1+\dim \Sl(\CC)=4$.
Moreover, it follows from
Lemma~\ref{lemma:notrivialonboundary}
and Remark~\ref{rem:nokleinonboundary} that the representation
$\rho^{+}$ verifies the hypothesis
of Lemma~\ref{lemma:obstructionsvanish}.
Hence all
cocycles in $Z^1(\p M; \Sl(\CC)_{\rho^+})$ are integrable. By
applying Artin's theorem \cite{Art68} we obtain from a formal
deformation
of $\rho^+$ a convergent deformation (see \cite[Lemma~3.3]{HPS01def}).
Thus $\rho^+$ is a smooth point of $R(M)$ with local
dimension equal to $4=\dim Z^1(\p M;{\Sl(\CC)}_{\rho^+})$.
\end{proof}

\section{The quadratic cone at the representation $\rho_\alpha$}
\label{section:quadraticcone}

Let $\alpha\co\p M \to \CC^{*}$ be a nontrivial representation. We
shall suppose in the sequel that $\alpha$ is a simple zero of the
Alexander invariant.

We want to show that $\rho_\alpha$ is contained in precisely two
components and that their intersection at the orbit  of
$\rho_\alpha$
is transverse. For this we study the quadratic cone. The Zariski
tangent space  can be viewed as the space of germs of
analytic paths which are contained in $R(M)$ at the  first order.
The quadratic  cone is the analogue at the second order. Since the
defining polynomials for the union of two varieties are the
products of defining polynomials for each component, the Zariski
tangent space of each component (at points of the intersection)
can only be detected by the quadratic cone.

Let $\rho\co\p M\to \PSL(\CC)$ be a representation.
The quadratic cone
$Q(\rho)$ is defined by
$$Q(\rho):=\{ u\in Z^{1}(\p M; \Sl(\CC)_{\rho})\mid [u\cup u]\sim
0\}\,.$$
 Recall that given two cocycles $u,v \in Z^1(\p M;
\Sl(\CC)_{\rho_\alpha})$
the cup product $[u\cup v]\in Z^2(\p M;{\Sl(\CC)}_{\rho_\alpha})$
is the cocycle  given by
\[
 [u\cup v](\gamma_1,\gamma_2)=
[u(\gamma_1),\Ad_{\rho_\alpha(\gamma_1)}(u(\gamma_2))]\in\Sl(\CC),
\qquad \forall\gamma_1\gamma_2\in\p M;
\]
where $[\,,\,]$ denotes the Lie bracket (see (\ref{eq:cupprod})).
Since the Lie bracket is antisymmetric, one easily checks that the
cup product is symmetric, i.e.\ the cocycles $[u\cup v]$ and
$[v\cup u]$ represent the same cohomology class, by
 (\ref{eq:cupsym}).

To compute the  quadratic cone  $Q(\rho_{\alpha})$ we use the
decomposition $\Sl(\CC)_{\alpha}=\CC_{0}\oplus \CC_+\oplus \CC_-$
of $\p M$-modules, see  (\ref{eqn:decomp}). Let
$\operatorname{pr}_0$ and $\operatorname{pr}_{\pm}$ denote the
projections of $\Sl(\CC)_{\alpha}$ to the respective one
dimensional modules. We can easily check the relation
\begin{equation}
\label{eqn:cupprojection}
\operatorname{pr}_{\pm}([u\cup v])=\operatorname{pr}_0(u)\overset{.}{\cup}
\operatorname{pr}_{\pm} (v)- \operatorname{pr}_{\pm} (u)\overset{.}{\cup}
\operatorname{pr}_0(v).
\end{equation}
Here we use the  products of $\p M$-modules $\CC\times \CC_{\pm}\to
\CC_{\pm}$, and $ \CC_{\pm}\times \CC \to \CC_{\pm}$, which are nothing but
the usual product of complex numbers. Notice that these cup
products of cohomology classes of cocycles valued in $\CC_{\pm}$ and
$\CC$ are antisymmetric (for the same reason that the cup product
valued in $\Sl(\CC)$ is symmetric).
 Therefore (\ref{eqn:cupprojection}) in
induces $\forall z\in  Z^1(\p M;{\Sl(\CC)}_{\alpha})$:
\begin{equation}
\label{eqn:cupprojectioncohomology} \operatorname{pr}_{\pm}([z\cup
z])=2\operatorname{pr}_0(z)\overset{.}{\cup} \operatorname{pr}_{\pm}
(z)\quad\textrm{up to coboundary}.
\end{equation}
The splitting (\ref{eqn:decomp}) induces a splitting in
cohomology. We recall that 
\[
\dim H^1(\p M;\CC_{\pm})=\dim H^1(\p M;\CC)=1
\]
and we have chosen cocycles $d_{\pm}$ and $d_0$ whose
cohomology classes generate $H^1(\p M;\CC_{\pm})$ and $ H^1(\p M;\CC)$
respectively (see Section~\ref{sec:cohomologymetabelian}).
Thus, the cocycle $z\in Z^1(\p M; \Sl(\CC)_{\alpha})$
is cohomologous to
\begin{equation}
\label{eqn:alphas}  z\sim a_0\,d_0+ a_-\,d_- + a_+\,d_+
\qquad\textrm{ with } a_0, a_+, a_-\in\CC.
\end{equation}
Therefore (\ref{eqn:cupprojectioncohomology}) becomes
\begin{equation}
\label{eqn:projectionalphas} \operatorname{pr}_{\pm}([z\cup z]) = 2
a_0 a_{ \pm}\, d_0\overset . \cup d_{\pm} \in Z^1(\p M;\CC_{\pm})\quad
\textrm{ up to coboundary.}
\end{equation}

\begin{lemma}
\label{rem:quadraticcone} The quadratic cone
$Q(\rho_{\alpha})\subset Z^1(\p M;{\Sl(\CC)}_{\alpha})$
is the union of two subspaces, one of dimension 4 and
another one of dimension 3. These subspaces  are precisely the
kernels of the projections
\[
\operatorname{pr}_{0}\co  Z^1(\p M;{\Sl(\CC)}_{\alpha}) \to 
Z^1(\p M;\CC) = H^1(\p M;\CC)
\]
and
\[
Z^1(\p M;{\Sl(\CC)}_{\alpha}) \stackrel {\scriptscriptstyle\operatorname{pr}_{+}\oplus
\operatorname{pr}_{-}}{\longrightarrow} Z^1(\p M;\CC_+\oplus\CC_-) 
\to H^1(\p M;\CC_+\oplus\CC_-)
\]
respectively. In particular, the intersection of these subspaces is the
space of coboundaries $B^{1}(\p M;\Sl(\CC)_{\alpha})$.
\end{lemma}

\begin{proof}
Notice that the space
\[Z^1(\p M;{\Sl(\CC)}_{\alpha}) =
Z^1(\p M;\CC_{0})\oplus Z^1(\p M;\CC_+\oplus\CC_-)\]
 is five
dimensional and that
$$B^{1}(\p M;{\Sl(\CC)}_{\alpha}) = B^{1}(\p M;\CC_{+}\oplus\CC_{-}).$$
Every cocycle $z\in Z^1(\p M;{\Sl(\CC)}_{\alpha})$ can be uniquely
written as $z = a_{0} d_{0} + a_{+} d_{+} + a_{-}
d_{-} +b$ where $b\in B^{1}(\p M;\CC_{+}\oplus\CC_{-}) $.
Combining Lemma~\ref{lemma:nonzeroproduct},
Equations~(\ref{eqn:alphas}) and (\ref{eqn:projectionalphas}), the
condition for the quadratic cone $[z\cup z]\sim 0$ reduces to:
\begin{equation*}
a_0\,a_+= a_0\,a_-=0\,.
\end{equation*}
This is of course a necessary condition for integrability. In
particular  we deduce that $z\in Q(\rho_{\alpha})$  if and only if
$z\in Z^{1}(\p M;\CC_{+}\oplus\CC_{-}) $ or
$z\in Z^1(\p M;\CC_{0})\oplus B^{1}(\p M;\CC_{+}\oplus\CC_{-})$
and the lemma follows.
\end{proof}

\begin{proof}[Proof of Theorem~\ref{thm:structureofR}] By
Lemma~\ref{rem:quadraticcone}  it suffices to show that $\rho_\alpha$
is contained in two
irreducible components, one of dimension four containing irreducible
representations and another of dimension three containing only
abelian ones. 
%As we said at the beginning of this section, the Zariski tangent space of
%each component of a union is only detected at the quadratic cone of
%the union (tangent to a point of the intersection). 
This will show
that the Zariski tangent space of each component has the right
dimension which implies that the point of the intersection is a
smooth point of each component. Moreover, the intersection is
transverse.

 The component of dimension 4 is provided by
Proposition~\ref{prop:rho+smooth}. In fact $\rho_\alpha$ is the adherence
of the orbit
of $\rho^+$, thus it is contained in the same irreducible component.
For the other component,
notice that $\rho_\alpha$ is
contained in a subvariety of abelian representations $S_{\alpha}(M)$, with $\dim S_{\alpha}(M) \geq 3$
(see Lemma~\ref{lem:nondef}).
Representations in $S_{\alpha}(M)$ are conjugate to diagonal representations,
thus the tangent
space to this deformation is clearly contained in the kernel of the
projection to $H^1(\p M;\CC_+\oplus\CC_-)$, and this gives an
irreducible component of dimension at most three.
Thus $\rho_{\alpha}$ it is a smooth point of the
three-dimensional component $S_{\alpha}(M) \subset R(M)$.
In addition,
the orbits by conjugation
in this component must be two-dimensional and therefore
all representations must be abelian.
\end{proof}

\begin{example} \label{ex:3iemepartie}
Let $M$ be the torus bundle given in Example~\ref{ex:1erepartie}.
Following Example~\ref{ex:2iemepartie} we choose
$\lambda_{i}\in\CC^{*}$ such that
 $\alpha_{i}\co\p M\to\CC^{*}$, $i=1,\ldots,4$, given by
$\alpha_{i}(\mu)=\lambda_{i}$ and $\alpha_{i}\big|_{\ZZ/2\oplus\ZZ/2}
 = \sigma_{i}$ is a simple root of the Alexander
invariant, i.e.\ $\lambda_{1}=3\pm \sqrt 8$ and $\lambda_{i}=1$ for
$i=2,3,4$. Therefore, in each case the representation
$\rho_{\alpha_{i}}$ is the limit of irreducible representations.
The deformation for $\rho_{\alpha_{4}}$ was already observed in
\cite[Section~4.2]{HP04}. This deformation was simply obtained as
a pullback of a component of the representation variety $R(\ZZ/2
\ast \ZZ/2)$ under the surjection:
\[\p M  
%=\pres{\alpha,\beta,\mu}{ \mu\alpha\mu^{-1}=\alpha\beta^{2},
%\mu\beta\mu^{-1}=\beta(\alpha\beta^{2})^{2}} 
\to
\p M / \langle \mu = 1 \rangle \cong \pres{\alpha,
\beta}{\alpha^{2},\beta^{2}}\,.
\]
Note that $\rho_{\alpha_{i}}$ is $\partial$-trivial for $i=2,3,4$.
On the other hand, $\rho^{+}_{\alpha_{i}}\circ i_{\#}\co
\p{\partial M}\to\PSL(\CC)$ is parabolic but nontrivial 
(cf.\ Lemma~\ref{lemma:notrivialonboundary}) and
Lemma~\ref{lemma:obstructionsvanish} applies. The results of
\cite{HPS01def} do not apply.
\end{example}

\section{The variety of characters near $\chi_\alpha$}
\label{section:characters}

Let $\alpha\co\p M\to \CC^{*}$ a representation such that
$\alpha$ is a simple zero of the Alexander invariant.
Let $\chi_\alpha\in X(M)$ denote the character of $\rho_\alpha$.

\begin{proof}[Proof of Theorem~\ref{thm:structureofX}]
Notice that there are at least two components of $X(M)$ containing
$\chi_\alpha$, which are precisely the quotients of the components
of $R(M)$ containing $\rho_\alpha$.

 To study the geometry of $X(M)$ near $\chi_\alpha$ we construct a slice as
in \cite{Lei02}. Let $\gamma_0\in\Gamma$ be an element such that
$\rho_\alpha(\gamma_0)\neq\pm I$. We define the slice as:
$$
{\mathcal S}=\{\rho\in R(M)\mid \rho(\gamma_0)\textrm{ is a
diagonal matrix }\}.
$$
By \cite{Lei02}, $\mathcal{S}$ is a slice \'etale and we
shall give here some of its properties:
firstly, ${\mathcal S}$ is transverse to the orbit by
conjugation at $\rho_\alpha$. More precisely,
if for some neighborhood $\mathcal U$ of $\rho_\alpha$, $F\co\mathcal U\subset R
(M)\to\CC^2$ maps a
representation $\rho$ to the non-diagonal entries of
$\rho(\gamma_0)$, then ${\mathcal S}\cap \mathcal U=F^{-1}(0,0)$. 
The tangent map of $F$
restricted to the coboundary space defines an isomorphism
$B^1(\p M;{\Sl(\CC)}_{\alpha})\cong \CC^2$. 
Thus
$$
 \TaZ{\rho_\alpha}{{\mathcal S}}\oplus B^1(\p M;{\Sl(\CC)}_{\alpha})=
Z^1(\p M;{\Sl(\CC)}_{\alpha}), $$ and $H^1(\p
M;{\Sl(\CC)}_{\alpha})$ can be viewed as the Zariski tangent
space to ${\mathcal S}$ at $\rho_\alpha$. Secondly, the projection
${\mathcal S}\to X(M)$ is surjective at least for characters
$\chi$ with $\chi(\gamma_0)\neq 4$. It is therefore sufficient for our purpose
to study the slice ${\mathcal S}$ and its quotient
by the stabilizer of $\rho_\alpha$.

Let $d_0\in Z^1(\p M;\CC_0)$ and $d_{\pm}\in Z^1(\p M;\CC_{\pm})$ denote
the cocycles of the previous section. Up to adding a coboundary,
we may assume that $d_{\pm}(\gamma_0)=0$. Thus the tangent space to
${\mathcal S}$ at $\rho_\alpha$ is three dimensional and generated
by $d_0$, $d_+$ and $d_-$. The analysis of
Section~\ref{section:quadraticcone} allows to say that ${\mathcal
S}$ has two components around $\rho_\alpha$: one curve tangent to
$d_0$ consisting of diagonal representation and a surface tangent
to the space generated by $d_+$ and $d_-$, containing irreducible
representations. Notice that since $\mathcal S$ is transverse to
the boundary space, it intersects the orbit of $\rho_{\alpha}$ by
conjugation in a single point.

To analyze $X(M)$, we take the quotient of ${\mathcal S}$ by the stabilizer of
$\rho_\alpha$, which is precisely the group of diagonal matrices: $$
D=\left\{\pm
\begin{pmatrix} \lambda & 0 \\ 0 & \lambda^{-1}
\end{pmatrix}\mid \lambda\in\CC^*\right\}\subset \PSL(\CC). $$
 Since the curve in ${\mathcal S}$ of diagonal representations  commutes with
$D$, it projects to a curve of abelian characters in
$X(M)$. To understand the action on the other component, we analyze
the action on the Zariski tangent space: a matrix
$\left(
\begin{smallmatrix} \lambda & 0 \\ 0 & \lambda^{-1}
\end{smallmatrix}\right)$ acts by mapping $d_{\pm}$ to
$\lambda^{\pm 2} d_{\pm}$. In other words, we must understand the
\emph{algebraic} quotient $\CC^2/\!/\CC^*$ where $t\in \CC^*$ maps
 $(x,y)\in
\CC^2$ to $(t \,x, t^{-1}\,y)$. The quotient $\CC^2/\!/\CC^*$ is
the line $\CC$, where the algebra of invariant functions is
generated by $xy$. Thus the quotient of the four dimensional
component containing $\rho_\alpha$ is also a smooth curve.

To show that the intersection is transverse, notice that  the Zariski
tangent space to $X(M)$ at $\chi_\alpha$ is
$H^1(\p M;\CC)\oplus
H^1(\p M;\CC_+\oplus \CC_-)/\!/\CC^*\cong \CC^2$. The first factor
$H^1(\p M;\CC)\oplus 0\cong \CC$ is tangent to the
curve of abelian characters, and $0\oplus H^1(\p M;\CC_+\oplus
\CC_-)/\!/\CC^*\cong \CC$ is tangent to the other curve, thus the
intersection is transverse.
\end{proof}

\section{Real valued characters}
\label{section:realvalued}

Let $\alpha\co\p M\to \CC^{*}$ be a representation such that
$\alpha$ is a simple zero of the Alexander invariant. Moreover, we
shall suppose in this section that the character $\chi_{\alpha}\co\p M \to\RR$ is
real valued.

We recall from \cite{HP04}
 that the character $\chi_\rho$ of a representation $\rho\in R(M)$ maps
$\gamma\in\pi_1(M)$ to
$\chi_\rho(\gamma)=\operatorname{trace}^2(\rho(\gamma))  $.

\begin{lemma}
Let $\Gamma$ be a finitely generated group. If $\chi\in X(\Gamma)$
is real valued, then there exists a representation $\rho\in
R(\Gamma)$ with character $\chi_{\rho}=\chi$ and such that  the
image of $\rho$ is contained either in $\mathrm{PSU}(2)$ or in
$\mathrm{PGL}_2(\RR)$.
\end{lemma}

Notice that
$\mathrm{PGL}_2(\RR)\subset\mathrm{PGL}_2(\CC)\cong\PSL(\CC)$
has two components, the
identity component is $\mathrm{PSL}_2(\RR)$,  the other component
consists of matrices with determinant $-1$ (which in
$\mathrm{PSL}_2(\CC)$ are represented
 by matrices with entries in $\CC$ with zero real part).

Looking at the action of $\PSL(\CC)$ on hyperbolic space $\HH^3$
by orientation preserving isometries, the group
$\mathrm{PGL}_2(\RR)$ is the stabilizer of a totally geodesic
plane in $\mathbb H^3$, and the two components of the group are
determined by whether their elements preserve or reverse 
the orientation of the
plane. The group $\mathrm{PSU}(2)$ is the stabilizer of a point
and it is connected (isomorphic to $\mathrm{SO}(3)$).

\begin{proof}
Let $F_n$ be a free group with a surjection
$F_n\twoheadrightarrow\Gamma$. Following \cite{HP04}, we consider
$F_n^{2}$, the subgroup of $F_n$  generated by all squared
elements $\gamma^2$, with $\gamma\in F_n$. This group  is a normal
subgroup of $F_{n}$ and gives the following exact sequence
\[
1\to F_n^{2}\to F_n\to H_1(F_n,C_2)\to 1,
\]
where $C_2$ is the cyclic group with two elements.

Let $\rho\in R(M)$ be a representation with real valued character
$\chi_{\rho}=\chi$. If $\rho$ is reducible, we may assume that it
is diagonal. Therefore we have two cases, either $\rho$ is
irreducible or $\rho$ is diagonal.

The composition of $\rho$ with $F_n\twoheadrightarrow\Gamma$ 
lifts to a representation $\rho'\co F_n\to \SL(\CC)$. By
Propositions~2.2 and 2.4 of \cite{HP04} the restriction of $\rho'$
to $F_n^2$ has real trace.  Hence, we can apply the known results
about these representations.

If the restriction $\rho'\vert_{F_n^2}$ is irreducible, then the
image of $\rho'\vert_{F_n^2}$ is contained in either
$\mathrm{SL}_2(\RR)$ or $\mathrm{SU}(2)$
\cite[Prop.~III.1.1]{MS84}. Looking at the action on $\mathbb
H^3$, this means that $\rho'\vert_{F_n^2}$ preserves either a
totally geodesic plane or a point in $\HH^3$. In particular
$\rho'\vert_{F_n^2}$ has either a unique invariant circle in
$\partial \HH^3\cong P^1 (\CC)$ or a unique fixed point in $\HH^3$
(uniqueness follows from irreducibility). Since $F^2_n$ is normal
in $F_n$, $\rho'(F_n)$ leaves invariant the same circle or the
same point, respectively. This proves the lemma in this case.

If the  restriction $\rho'\vert_{F_n^2}$ is trivial, then the
image of $\rho$ is abelian and finite, and therefore it fixes a
point in $\HH^3$. This means that, up to conjugation, $\rho(\Gamma)\subset
\mathrm{PSU}(2)$.

If the restriction $\rho'\vert_{F_n^2}$ is reducible but
non-trivial, then it fixes either a single point or two in
$\partial \HH^3$. When it fixes two points, using normality, these
points must be fixed or permuted by every element of
 $\rho'(F_n)$, and therefore
$\rho$ is either diagonal or the Klein group and the lemma is
shown in this case. Finally, if $\rho'\vert_{F_n^2}$ fixes a
single point, it is also fixed by $\rho$, but this contradicts the
fact that $\rho$ is diagonal.
\end{proof}

As $\chi_{\alpha}$ is real valued,
$\alpha(\gamma)+\frac1{\alpha(\gamma)}\in\RR$
$\forall\gamma\in\pi_1(M)$.
Thus the image of $\alpha$ is either contained in the reals $\RR^*$ or
in the unit circle $S^1=\{z\in\CC\mid \vert z\vert=1\}$.

 Since we already
know which cocycles of $Z^{1}(\p M;\Sl(\CC)_{\alpha})$ are
integrable, we can easily describe the subsets
 of representations into these groups. We distinguish two cases:

\medskip

\textbf{Case 1}\qua Assume that there is $\gamma\in\p M$ such that
$|\alpha(\gamma)|\neq 1$.

\medskip

In this case the image of $\alpha$ is contained in $\RR^*=\RR\setminus
\{0\}$ but not in $\{\pm 1\}$. Hence, $\image(\rho_\alpha)  $ is
contained
in $\mathrm{PGL}_2(\RR)$. Since $\Sl(\CC)$ is the complexification of
$\Sl(\RR)$, we have the corresponding isomorphism of cohomology groups:
$$
H^*(\p M;\Sl(\CC)_{\alpha})=
H^*(\p M;{\Sl(\RR)}_{\alpha
})\otimes_{\RR}\CC. $$
In particular, we may assume that the cocycles $d_{\pm}$ and $d_0$ are
valued in $\mathbb R$. Using the complexification of
the second cohomology group, we realize that the computation of
obstructions can be carried out in the real setting, thus we
get:

\begin{prop}
\label{prop:realPSL}
Assume that the image of $\alpha$ is contained in
$\RR^*$ but not in $\{\pm 1\}$. Then the character $\chi_\alpha$ is contained in
precisely two real curves of characters in
$\mathrm{PGL}_2(\RR)$, one of them with abelian representations and
the other one with irreducible ones. In addition
$\chi_\alpha$ is a smooth point of both curves that intersect transversely
at  $\chi_\alpha$.
\end{prop}

When the image of $\alpha$ is  contained in the positive reals,
then the representations in Proposition~\ref{prop:realPSL} are
contained in $\PSL(\RR)$, by connectedness. The case when the
image of $\alpha$ is contained in $\{\pm 1\}$, is treated in next
case.

\medskip

\textbf{Case 2}\qua Assume that the image of $\alpha$ is
contained in $S^1=\{z\in\CC\mid \vert z\vert=1\}$.

\medskip

Now $\image(\rho_\alpha)$ is contained in the intersection
$\mathrm{PSU}(2)\cap\mathrm{PSU}(1,1)$, where
$$
\mathrm{PSU}(2)=\pm\left\{
\begin{pmatrix} a & b\\ -\overline b & \overline a
\end{pmatrix}\mid \vert a\vert^2+\vert
b\vert^2=1  \right\}
$$ and $$
\mathrm{PSU}(1,1)=\pm\left\{
\begin{pmatrix} a & b\\ \overline b & \overline a
\end{pmatrix}\mid \vert a\vert^2-\vert
b\vert^2=1  \right\}. $$
Geometrically, $\mathrm{PSU}(2)$ is the stabilizer of a point in
hyperbolic space and $\mathrm{PSU}(1,1)$ is the connected
component of the stabilizer of the unit circle in $\CC$. Thus
$\mathrm{PSU}(1,1)$ is the connected component of the
stabilizer of a plane in hyperbolic space, and therefore it is
conjugate to $\PSL(\RR)$.

In this case, the Lie algebra $\Sl(\CC) $ is the complexification of
both $\su$ and $\mathfrak{su}(1,1)$. Thus an argument
similar to
the previous case will apply. Note that the cocycles $d_+$ and $d_-$
are not valuated in those Lie algebras. However, since 
$\ol{\rho_\alpha} = \rho_{1/\alpha}$ we can assume
that  $d_-$ is the complex conjugate transpose to $d_+$, thus
$d_+-d_-$ is valuated in $\su$ and
$d_++d_-$ is valuated in $\mathfrak{su}(1,1)$.

The tangent directions  $d_+-d_-$ and  $d_++d_-$ are different in the
slice of the variety of representations, but they
project to the same direction (with opposite sense) in the variety of
characters (see the description of the quotient in the
previous section). Notice that this gives a curve of real valued
characters. In addition, the set of real valued characters in a smooth
complex curve can be  at most one dimensional.

\begin{prop}
Assume that the image of $\alpha$ is contained in
$S^1\subset\CC^*$. Then the character $\chi_\alpha$ is contained
in precisely two real curves of characters, one of them abelian
(contained in $S^1$) and the other one with irreducible
representations, contained in $\mathrm{PSU}(2)$ in one side and in
$\mathrm{PSU}(1,1)\cong\PSL(\RR)$ on the other side.
 In addition, $\chi_\alpha$ is a
smooth point of both curves that intersect transversely at
$\chi_\alpha$.
\end{prop}

\Addresses\recd

\end{document}